\documentclass[11pt,twoside]{article}
\usepackage{theorem}
\usepackage{epsfig,amsmath,amssymb}
\usepackage{color}
\allowdisplaybreaks[1]

\newtheorem{thm}{\bf Theorem}[section]

\newtheorem{lem}[thm]{\bf Lemma}

\newtheorem{defn}{\bf Definition}[section]
\newtheorem{rem}{\bf Remark}[section]
\numberwithin{equation}{section}

\newtheorem{prob}{\bf Problem}[section]

\title{\LARGE\bf Extension of isometries in real Hilbert spaces}

\author{Soon-Mo Jung}

\date{\small\it Mathematics Section,
                College of Science and Technology, \\
                Hongik University, 30016 Sejong, Korea \\
                E-mail: {\tt smjung@hongik.ac.kr} \\
                \hspace*{12.5mm} {\tt smjung747@gmail.com}}

\begin{document}

\markboth{Extension of isometries}{S.-M. Jung}
\maketitle

\begin{abstract}
In this paper, the notions of first-order and second-order
generalized linear spans and index set are defined.
Moreover, their properties are investigated and applied to the
studies of extension of isometries.
We develop the theory of extending the domain of local
isometries to the generalized linear spans, where we call an
isometry defined in a subset of a Hilbert space a local isometry.
In addition, we prove that the domain of local isometry can be
extended to any real Hilbert space, where the domain of local
isometry does not have to be a convex body or an open set.
This indicates that the main results of this paper are superior
to those previously published.
\end{abstract}
\vspace{5mm}

\noindent
{\bf AMS Subject Classification:} Primary 46B04; Secondary 46C99.
\vspace{5mm}

\noindent
{\bf Key Words:} isometry; extension of isometry;
                 generalized linear span; index; index set.
\normalsize


\section{Introduction}
\label{sec:1}

In the course of the development of mathematics in the last
century, the problem of extending the domain of a function
while maintaining the characteristics of a function defined in
a local domain had a particularly great influence on the
development of functional analysis.

For example, in topology, the Tietze extension theorem states
that all continuous functions defined on a closed subset of a
normal topological space can be extended to the entire space.

\begin{thm}[Tietze] \label{thm:tietze}
Let $X$ be a normal space, $E$ a nonempty closed subset of $X$,
and let $[-L, L]$ be a closed real interval.
If $f : E \to [-L, L]$ is a continuous function, then there
exists a continuous extension of $f$ to $X$, \textit{i.e.},
there exists a continuous function $F : X \to [-L, L]$ such
that $F(x) = f(x)$ for all $x \in E$.
\end{thm}

The Tietze extension theorem has a wide range of applications
and is an interesting theorem, so there are many variations of
this theorem.

As another example, the following theorem is widely known
regarding the extension of continuous linear functions.

\begin{thm} \label{thm:clm}
Let $X$ and $Y$ be a normed space and a Banach space,
respectively.
If $E$ is a subspace of $X$, then every continuous linear
function $f : E \to Y$ has a unique extension to a continuous
linear function $F : \overline{E} \to Y$, where $\overline{E}$
is the closure of $E$ in $X$.
\end{thm}

In 1972, Mankiewicz published his paper \cite{manki} to
determining whether an isometry $f : E \to Y$ from a subset
$E$ of a real normed space $X$ into a real normed space $Y$
admits an extension to an isometry from $X$ onto $Y$.
Indeed, he proved that every surjective isometry $f : E \to Y$
can be uniquely extended to an affine isometry between the
whole spaces when either $E$ and $f(E)$ are both convex bodies
or $E$ is nonempty open connected and $f(E)$ is open, where
convex body is a convex set with nonempty interior.

\begin{thm}[Mankiewicz] \label{thm:mankiewicz}
Let $X$ and $Y$ be real normed spaces, $E$ a nonempty subset
of $X$, and let $f : E \to f(E)$ be a surjective isometry,
where $f(E)$ is a subset of $Y$.
If either both $E$ and $f(E)$ are convex bodies or $E$ is open
and connected and $f(E)$ is open, then $f$ can be uniquely
extended to an affine isometry $F : X \to Y$.
\end{thm}

This conclusion particularly holds for the closed unit balls.
Based on this fact, with the same research direction, Tingley
\cite{tingley} intuitively paid attention to the unit spheres
and posed the following problem, which is now known as
Tingley's problem.

\begin{prob}[Tingley] \label{prob:tingley}
Is every surjective isometry between the unit spheres of two
Banach spaces the restriction to the unit sphere of a
surjective linear isometry between the whole spaces?
\end{prob}

Recently, papers on the extension of isometries related to
Tingley's problem have been mainly published
(see \cite{ding2002,ding2003,ding2009,ding2015,vershik}).

The research in this paper is strongly motivated by Theorem
\ref{thm:mankiewicz} and Problem \ref{prob:tingley}, among
others (refer to \cite{avellan,guang,solecki} also).

In Section \ref{sec:3} of this paper, we introduce new concepts
such as first-order generalized linear span and index set,
which are essential to prove the final result of this paper.
Section \ref{sec:4} is devoted to developing the theory of
extending the domain of local isometries to the first-order
generalized linear span, where we call an isometry defined in
a subset of a Hilbert space a local isometry.
We introduce the concept of second-order generalized linear
span in Section \ref{sec:5} and develop the theory of extending
the domain of local isometries to the second-order generalized
linear span in Section \ref{sec:7}.
Finally, we prove in Theorem \ref{thm:3.17} that the domain of
local isometry can be extended to any real Hilbert space
including that domain, where the domain of a local isometry
needs to be bounded and have at least two elements, but need
not be a convex body or an open set.
This indicates that the main results of this paper are superior
to those previously published.


\section{Preliminaries}
\label{sec:2}

Throughout this paper, using the symbol $\mathbb{R}^\omega$,
we represents an infinite dimensional real vector space
defined as
\begin{eqnarray*}
\mathbb{R}^\omega = \big\{ (x_1, x_2, \ldots ) :
x_i \in \mathbb{R} ~\,\mbox{for all}\,~ i \in \mathbb{N} \,\big\}.
\end{eqnarray*}
From now on, we denote by $(\mathbb{R}^\omega, \mathcal{T})$ the
product space $\prod\limits_{i=1}^\infty \mathbb{R}$, where
$(\mathbb{R}, \mathcal{T}_{\mathbb{R}})$ is the usual topological
space.
Then since $(\mathbb{R}, \mathcal{T}_{\mathbb{R}})$ is a
Hausdorff space, $(\mathbb{R}^\omega, \mathcal{T})$ is a
Hausdorff space (see \cite[Theorem 111.7]{kasriel}).

Let $a = \{ a_i \}$ be a sequence of positive real numbers
satisfying the condition
\begin{eqnarray} \label{eq:20110922-1}
\sum_{i=1}^\infty a_i^2 < \infty.
\end{eqnarray}
With this sequence $a = \{ a_i \}$, we define
\begin{eqnarray*}
M_a = \Bigg\{ (x_1, x_2, \ldots) \in \mathbb{R}^\omega :\,
              \sum_{i=1}^\infty a_i^2 x_i^2 < \infty
      \Bigg\}.
\end{eqnarray*}
Then $M_a$ is a vector space over $\mathbb{R}$.
And we can define an inner product
$\langle \cdot, \cdot \rangle_a$ for $M_a$ by
\begin{eqnarray*}
\langle x, y \rangle_a = \sum_{i=1}^\infty a_i^2 x_i y_i
\end{eqnarray*}
for all $x = (x_1, x_2, \ldots)$ and $y = (y_1, y_2, \ldots)$
of $M_a$, with which $(M_a, \langle \cdot, \cdot \rangle_a)$
becomes a real inner product space.
This inner product induces the norm in the natural way
\begin{eqnarray*}
\| x \|_a = \sqrt{\langle x, x \rangle_a}
\end{eqnarray*}
for all $x \in M_a$, so that $(M_a, \| \cdot \|_a)$ becomes a
real normed space.


\begin{rem} \label{rem:new2.1}
$M_a$ is the set of each element $x \in \mathbb{R}^\omega$
that satisfies $\| x \|_a^2 < \infty$, \textit{i.e.},
\begin{eqnarray*}
M_a = \big\{ (x_1, x_2, \ldots) \in \mathbb{R}^\omega :\,
             \| x \|_a^2 < \infty
      \big\}.
\end{eqnarray*}
\end{rem}

Furthermore, we define the metric $d_a$ on $M_a$ by
\begin{eqnarray*}
d_a(x,y) = \| x - y \|_a = \sqrt{\langle x-y, x-y \rangle_a}
\end{eqnarray*}
for all $x, y \in M_a$.
Thus, $(M_a, d_a)$ is a real metric space.

What follows is a basic definition we are familiar with, but
for the sake of completeness of the paper, we now define
precisely the $d_a$-isometry between subsets of $M_a$.


\begin{defn} \label{def:1.1}
Let $E_1$ and $E_2$ be nonempty subsets of $M_a$.
\begin{itemize}
\item[$(i)$]  A function $f : E_1 \to E_2$ is called a
              $d_a$-\emph{isometry} provided
              $d_a(f(x), f(y)) = d_a(x,y)$ for all
              $x, y \in E_1$;
\item[$(ii)$] $E_1$ is said to be $d_a$-\emph{isometric} to
              $E_2$ provided there exists a surjective
              $d_a$-isometry $f : E_1 \to E_2$.
\end{itemize}
\end{defn}

Assuming that the sequence $a = \{ a_i \}$ satisfies the
condition (\ref{eq:20110922-1}), we can prove that
$(M_a, \langle \cdot, \cdot \rangle_a)$ is a Hilbert space in
the same way as \cite[Theorem 2.1]{jungkim}.
Let $(M_a, \mathcal{T}_a)$ be the topological space generated by
the metric $d_a$.


\begin{rem} \label{rem:2.1}
We note that
\begin{itemize}
\item[$(i)$]   $(M_a, \langle \cdot, \cdot \rangle_a)$ is a
               Hilbert space over $\mathbb{R}$;
\item[$(ii)$]  $(M_a, \mathcal{T}_a)$ is a Hausdorff space as
               a subspace of the Hausdorff space
               $(\mathbb{R}^\omega, \mathcal{T})$.
\end{itemize}
\end{rem}


\begin{defn} \label{def:2.3}
Given a $c \in M_a$, a \emph{translation} $T_c : M_a \to M_a$
is a function defined by $T_c (x) = x + c$ for all $x \in M_a$.
\end{defn}


\section{First-order generalized linear span}
\label{sec:3}

In \cite[Theorem 2.5]{jungkim}, we were able to extend the
domain of $d_a$-isometry $f$ to the entire space when the
domain of $f$ is a \emph{non-degenerate basic cylinder}
(see Definition \ref{def:basic} for the exact definition of
non-degenerate basic cylinders).
However, we shall see in Definition \ref{def:ext} and Theorem
\ref{thm:new3} that the domain of $d_a$-isometry $f$ can be
extended to its first-order generalized linear span whenever
$f$ is defined on a bounded set which contains more than one
element.

From now on, it is assumed that $E$, $E_1$, and $E_2$ are
subsets of $M_a$, each of them contains more than one element,
unless specifically stated for their cardinalities.
If the set has only one element or no element, this case will
not be covered here because the results derived from this case
are trivial and uninteresting.


\begin{defn} \label{def:gls}
Assume that $E$ is a nonempty bounded subset of $M_a$ and $p$
is a fixed element of $E$.
We define the \emph{first-order generalized linear span} of
$E$ with respect to $p$ as
\begin{eqnarray*}
\begin{split}
{\rm GS}(E,p)
& = \Bigg\{ p + \sum_{i=1}^m \sum_{j=1}^\infty
            \alpha_{ij} (x_{ij} - p) \in M_a :
            m \in \mathbb{N}\,;
            \vspace*{1mm}\\
& \hspace*{10mm}
            x_{ij} \in E ~\mbox{and}~
            \alpha_{ij} \in \mathbb{R}~
            \mbox{for~all~$i$~and~$j$}
    \Bigg\}.
\end{split}
\end{eqnarray*}
\end{defn}

We remark that if a bounded subset $E$ of $M_a$ contains more
than one element, then $E$ is a proper subset of its first-order
generalized linear span ${\rm GS}(E,p)$, because
$x = p + (x - p) \in {\rm GS}(E,p)$ for any $x \in E$ and
$p + \alpha (x - p) \in {\rm GS}(E,p)$ for any
$\alpha \in \mathbb{R}$, which implies that ${\rm GS}(E,p)$ is
unbounded.
Moreover, we note that $\alpha x + \beta y \in M_a$ for all
$x, y \in M_a$ and $\alpha, \beta \in \mathbb{R}$, because
$\| \alpha x + \beta y \|_a
 \leq | \alpha | \| x \|_a + | \beta | \| y \|_a < \infty$.
Therefore, ${\rm GS}(E,p) - p$ is a real vector space, because
the double sum in the definition of ${\rm GS}(E,p)$ guarantees
$\alpha x + \beta y \in {\rm GS}(E,p) - p$ for all
$x, y \in {\rm GS}(E,p) - p$ and $\alpha, \beta \in \mathbb{R}$
and because ${\rm GS}(E,p) - p$ is a subset of a real vector
space $M_a$ (\textit{cf.} Lemma \ref{lem:3.8} $(i)$ below).

For each $i \in \mathbb{N}$, we set
$e_i = (0, \ldots, 0, 1, 0, \ldots)$, where $1$ is in the
$i$th position.
Then $\big\{ \frac{1}{a_i} e_i \big\}$ is a complete
orthonormal sequence in $M_a$.


\begin{defn} \label{def:Lambda}
Given a nonempty subset $E$ of $M_a$, we define the
\emph{index set} of $E$ by
\begin{eqnarray*}
\Lambda(E) = \Big\{ i \in \mathbb{N} : \mbox{there are an}~
                    x \in E ~\mbox{and an}~
                    \alpha \in \mathbb{R} \!\setminus\! \{ 0 \}
                    ~\mbox{satisfying}~ x + \alpha e_i \in E
             \Big\}.
\end{eqnarray*}
Each $i \in \Lambda(E)$ is called an \emph{index} of $E$.
If $\Lambda(E) \neq \mathbb{N}$, then the set $E$ is called
\emph{degenerate}.
Otherwise, $E$ is called \emph{non-degenerate}.
\end{defn}

We will find that the concept of index set in Hilbert space
sometimes takes over the role that the concept of dimension
plays in vector space.
According to the definition above, if $i$ is an index of $E$,
\textit{i.e.}, $i \in \Lambda(E)$, then there are $x \in E$
and $x + \alpha e_i \in E$ for some $\alpha \neq 0$.
Since $x \neq x + \alpha e_i$, we remark that if
$\Lambda(E) \neq \emptyset$, then the set $E$ contains at least
two elements.

In the following lemma, we prove that if $i$ is an index of
$E$ and $p \in E$, then the first-order generalized linear span
${\rm GS}(E,p)$ contains the line through $p$ in the direction
$e_i$.


\begin{lem} \label{lem:2018-9-18}
Assume that $E$ is a bounded subset of $M_a$ and
${\rm GS}(E,p)$ is the first-order generalized linear span of
$E$ with respect to a fixed element $p \in E$.
If $i \in \Lambda(E)$, then $p + \alpha e_i \in {\rm GS}(E,p)$
for all $\alpha \in \mathbb{R}$.
\end{lem}

\noindent
\emph{Proof}.
By Definition \ref{def:Lambda}, if $i \in \Lambda(E)$ then there
exist an $x \in E$ and an $\alpha_0 \neq 0$, which satisfy
$x + \alpha_0 e_i \in E$.
Since $x \in E$ and $x + \alpha_0 e_i \in E$, by Definition
\ref{def:gls}, we get
\begin{eqnarray*}
p + \alpha_0 \beta e_i
=   p + \beta (x + \alpha_0 e_i - p) - \beta (x - p)
\in {\rm GS}(E,p)
\end{eqnarray*}
for all $\beta \in \mathbb{R}$.
Setting $\alpha = \alpha_0 \beta$ in the above relation, we
obtain $p + \alpha e_i \in {\rm GS}(E,p)$ for any
$\alpha \in \mathbb{R}$.
\hfill$\Box$
\vspace*{5mm}

In the following lemma, we introduce a characteristic property
of the index of a nonempty bounded subset of $M_a$.
We prove that $i$ is an index of $E$ if and only if
$x_i \neq p_i$ for some
$x = (x_1, x_2, \ldots) \in {\rm GS}(E,p)$.


\begin{lem} \label{lem:doolset}
Let $E$ be a bounded subset of $M_a$ and let
$p = (p_1, p_2, \ldots)$ be an element of $E$.
Then $i \not\in \Lambda(E)$ if and only if $x_i = p_i$ for all
$x = (x_1, x_2, \ldots) \in {\rm GS}(E,p)$.
\end{lem}

\noindent
\emph{Proof}.
First, we prove that $i \not\in \Lambda(E)$ if and only if
$x_i = p_i$ for all $x \in E$.
Assume that $x_i = p_i$ for all $x \in E$.
Then for any $y \in E$ and $\alpha \neq 0$, we have
$y + \alpha e_i \not\in E$ because the $i$th coordinate of
$y + \alpha e_i$ is $y_i + \alpha = p_i + \alpha \neq p_i$.
This fact implies that $i \not\in \Lambda(E)$.

Now we assume that $i \not\in \Lambda(E)$ and, on the contrary,
assume that $x_i \neq p_i$ for some $x \in E$.
Because $p$ and $x$ are elements of $E$, we have $p_j \in E_j$
and $x_j \in E_j$ for all $j \in \mathbb{N}$, where we set
$E_j := \{ y_j : (y_1, y_2, \ldots) \in E \}$.
Thus, we have
$p + (x_i - p_i) e_i
 = (p_1, \ldots, p_{i-1}, x_i, p_{i+1}, \ldots) \in E$ because
each $j$th coordinate of $p + (x_i - p_i) e_i$ belongs to $E_j$,
which would imply that $i \in \Lambda(E)$, that is a
contradiction.
Hence, we proved that if $i \not\in \Lambda(E)$ then $x_i = p_i$
for all $x \in E$.

Finally, in view of Definition \ref{def:gls}, it is easy to
see that $x_i = p_i$ for all $x \in E$ if and only if
$x_i = p_i$ for all $x \in {\rm GS}(E,p)$, which completes our
proof.
\hfill$\Box$
\vspace*{5mm}

We now introduce a lemma, which is a generalized version of
\cite[Lemma 2.3]{jungkim} and whose proof runs in the same way.
We prove that the function
$T_{-q} \circ f \circ T_p : E_1 - p \to E_2 - q$ preserves the
inner product.
This property is important in proving the following theorems
as a necessary condition for $f$ to be a $d_a$-isometry.


\begin{lem} \label{lem:new2}
Assume that $E_1$ and $E_2$ are bounded subsets of $M_a$
that are $d_a$-isometric to each other via a surjective
$d_a$-isometry $f : E_1 \to E_2$.
Assume that $p$ is an element of $E_1$ and $q$ is an element
of $E_2$ with $q = f(p)$.
Then the function $T_{-q} \circ f \circ T_p : E_1 - p \to E_2 - q$
preserves the inner product, \textit{i.e.},
\begin{eqnarray*}
\big\langle ( T_{-q} \circ f \circ T_p )(x-p),\,
            ( T_{-q} \circ f \circ T_p )(y-p)
\big\rangle_a = \langle x-p, y-p \rangle_a
\end{eqnarray*}
for all $x, y \in E_1$.
\end{lem}

\noindent
\emph{Proof}.
Since $T_{-q} \circ f \circ T_p : E_1 - p \to E_2 - q$ is a
$d_a$-isometry, we have
\begin{eqnarray*}
\| ( T_{-q} \circ f \circ T_p )(x-p) -
   ( T_{-q} \circ f \circ T_p )(y-p)
\|_a = \| (x-p) - (y-p) \|_a
\end{eqnarray*}
for any $x, y \in E_1$.
If we put $y = p$ in the last equality, then we get
\begin{eqnarray*}
\| ( T_{-q} \circ f \circ T_p )(x-p) \|_a = \| x-p \|_a
\end{eqnarray*}
for each $x \in E_1$.
Moreover, it follows from the previous equality that
\begin{eqnarray*}
\begin{split}
& \| ( T_{-q} \circ f \circ T_p )(x-p) -
     ( T_{-q} \circ f \circ T_p )(y-p)
  \|_a^2 \\
& \hspace*{3mm}
  = \big\langle ( T_{-q} \circ f \circ T_p )(x-p) -
                ( T_{-q} \circ f \circ T_p )(y-p), \\
& \hspace*{10mm}
                ( T_{-q} \circ f \circ T_p )(x-p) -
                ( T_{-q} \circ f \circ T_p )(y-p)
    \big\rangle_a \\
& \hspace*{3mm}
  = \| x-p \|_a^2 -
    2 \big\langle ( T_{-q} \circ f \circ T_p )(x-p),\,
                  ( T_{-q} \circ f \circ T_p )(y-p)
      \big\rangle_a + \| y-p \|_a^2
\end{split}
\end{eqnarray*}
and
\begin{eqnarray*}
\begin{split}
\| (x-p) - (y-p) \|_a^2
& = \big\langle (x-p) - (y-p), (x-p) - (y-p) \big\rangle_a \\
& = \| x-p \|_a^2 - 2 \langle x-p, y-p \rangle_a +
    \| y-p \|_a^2.
\end{split}
\end{eqnarray*}
Finally, comparing the last two equalities yields the validity
of our assertion.
\hfill$\Box$


\section{First-order extension of isometries}
\label{sec:4}

In the previous section, we made all the necessary preparations
to extend the domain $E_1$ of the surjective $d_a$-isometry
$f : E_1 \to E_2$ to its first-order generalized linear span
${\rm GS}(E_1,p)$.

Although $E_1$ is a bounded set, ${\rm GS}(E_1,p) - p$ is a
real vector space.
Now we will extend the $d_a$-isometry $T_{-q} \circ f \circ T_p$
defined on the bounded set $E_1 - p$ to the $d_a$-isometry
$T_{-q} \circ F \circ T_p$ defined on the vector space
${\rm GS}(E_1,p) - p$.
Comparing their `sizes' of $E_1 - p$ and ${\rm GS}(E_1,p) - p$,
or considering that ${\rm GS}(E_1,p) - p$ is an algebraically
closed space, it is a great achievement to extend the
$d_a$-isometry $T_{-q} \circ f \circ T_p$ defined on the
bounded set $E_1 - p$ to the $d_a$-isometry defined on the
vector space ${\rm GS}(E_1,p) - p$.


\begin{defn} \label{def:ext}
Assume that $E_1$ and $E_2$ are nonempty bounded subsets of
$M_a$ that are $d_a$-isometric to each other via a surjective
$d_a$-isometry $f : E_1 \to E_2$.
Let $p$ be a fixed element of $E_1$ and let $q$ be an element
of $E_2$ that satisfies $q = f(p)$.
We define a function $F : {\rm GS}(E_1,p) \to M_a$ as
\begin{eqnarray*} 
(T_{-q} \circ F \circ T_p)
\Bigg( \sum_{i=1}^m \sum_{j=1}^\infty
       \alpha_{ij} (x_{ij} - p) \Bigg)
= \sum_{i=1}^m \sum_{j=1}^\infty
  \alpha_{ij} (T_{-q} \circ f \circ T_p)(x_{ij} - p)
\end{eqnarray*}
for any $m \in \mathbb{N}$, $x_{ij} \in E_1$, and for all
$\alpha_{ij} \in \mathbb{R}$ satisfying
$\sum\limits_{i=1}^m \sum\limits_{j=1}^\infty
 \alpha_{ij} (x_{ij} - p) \in M_a$.
\end{defn}

We note that in the definition above, it is important for the
argument of $T_{-q} \circ F \circ T_p$ to belong to $M_a$.
Now we show that the function $F : {\rm GS}(E_1,p) \to M_a$ is
well defined.


\begin{lem} \label{lem:3.2}
Assume that $E_1$ and $E_2$ are bounded subsets of $M_a$
that are $d_a$-isometric to each other via a surjective $d_a$-isometry
$f : E_1 \to E_2$.
Let $p$ be an element of $E_1$ and let $q$ be an element of
$E_2$ that satisfy $q = f(p)$.
The function $F : {\rm GS}(E_1,p) \to M_a$ given in Definition
\ref{def:ext} is well defined.
\end{lem}

\noindent
\emph{Proof}.
First, we will check that the range of $F$ is a subset of $M_a$.
For any $m, n_1, n_2 \in \mathbb{N}$ with $n_2 > n_1$,
$x_{ij} \in E_1$, and for all $\alpha_{ij} \in \mathbb{R}$, it
follows from Lemma \ref{lem:new2} that
\begin{eqnarray} \label{eq:2020-1-10}
\begin{split}
& \bigg\| \sum_{i=1}^m \sum_{j=1}^{n_2}
          \alpha_{ij} (T_{-q} \circ f \circ T_p)(x_{ij} - p) -
          \sum_{i=1}^m \sum_{j=1}^{n_1}
          \alpha_{ij} (T_{-q} \circ f \circ T_p)(x_{ij} - p)
  \bigg\|_a^2 \\
& \hspace*{3mm} =
  \bigg\langle
  \sum_{i=1}^m \sum_{j = n_1 + 1}^{n_2}
  \alpha_{ij} (T_{-q} \circ f \circ T_p)(x_{ij} - p),\,
  \sum_{k=1}^m \sum_{\ell = n_1 + 1}^{n_2}
  \alpha_{k\ell} (T_{-q} \circ f \circ T_p)(x_{k\ell} - p)
  \bigg\rangle_a \\
& \hspace*{3mm} =
  \sum_{i=1}^m \sum_{k=1}^m
  \sum_{j = n_1 + 1}^{n_2} \alpha_{ij}
  \sum_{\ell = n_1 + 1}^{n_2} \alpha_{k\ell}
  \big\langle (T_{-q} \circ f \circ T_p)(x_{ij} - p),\,
              (T_{-q} \circ f \circ T_p)(x_{k\ell} - p)
  \big\rangle_a \\
& \hspace*{3mm} =
  \sum_{i=1}^m \sum_{k=1}^m
  \sum_{j = n_1 + 1}^{n_2} \alpha_{ij}
  \sum_{\ell = n_1 + 1}^{n_2} \alpha_{k\ell}
  \big\langle x_{ij} - p,\, x_{k\ell} - p \big\rangle_a \\
& \hspace*{3mm} =
  \bigg\langle
  \sum_{i=1}^m \sum_{j = n_1 + 1}^{n_2} \alpha_{ij}
  (x_{ij} - p),\,
  \sum_{k=1}^m \sum_{\ell = n_1 + 1}^{n_2} \alpha_{k\ell}
  (x_{k\ell} - p)
  \bigg\rangle_a \\
& \hspace*{3mm} =
  \bigg\|
  \sum_{i=1}^m \sum_{j = n_1 + 1}^{n_2} \alpha_{ij} (x_{ij} - p)
  \bigg\|_a^2 \\
& \hspace*{3mm} =
  \bigg\| \sum_{i=1}^m \sum_{j=1}^{n_2}
          \alpha_{ij} (x_{ij} - p) -
          \sum_{i=1}^m \sum_{j=1}^{n_1}
          \alpha_{ij} (x_{ij} - p)
  \bigg\|_a^2.
\end{split}
\end{eqnarray}
Indeed, the equality (\ref{eq:2020-1-10}) holds for all
$m, n_1, n_2 \in \mathbb{N}$.

We now assume that
$\sum\limits_{i=1}^m \sum\limits_{j=1}^\infty
 \alpha_{ij} (x_{ij} - p) \in M_a$ for some $x_{ij} \in E_1$
and $\alpha_{ij} \in \mathbb{R}$, where $m$ is a fixed positive
integer.
Then since $(M_a, \mathcal{T}_a)$ is a Hausdorff space on
account of Remark \ref{rem:2.1} $(ii)$ and the topology
$\mathcal{T}_a$ is consistent with the metric $d_a$ and with
the norm $\| \cdot \|_a$, the sequence
$\big\{ \sum\limits_{i=1}^m \sum\limits_{j=1}^n
        \alpha_{ij} (x_{ij} - p) \big\}_n$ converges to
$\sum\limits_{i=1}^m \sum\limits_{j=1}^\infty
 \alpha_{ij} (x_{ij} - p)$ (in $M_a$) and hence, the sequence
$\big\{ \sum\limits_{i=1}^m \sum\limits_{j=1}^n
        \alpha_{ij} (x_{ij} - p) \big\}_n$ is a Cauchy sequence
in $M_a$.

We know by (\ref{eq:2020-1-10}) and the definition of Cauchy
sequences that for each $\varepsilon > 0$ there exists an
integer $N_\varepsilon > 0$ such that
\begin{eqnarray*}
\begin{split}
& \bigg\| \sum_{i=1}^m \sum_{j=1}^{n_2}
          \alpha_{ij} (T_{-q} \circ f \circ T_p)(x_{ij} - p) -
          \sum_{i=1}^m \sum_{j=1}^{n_1}
          \alpha_{ij} (T_{-q} \circ f \circ T_p)(x_{ij} - p)
  \bigg\|_a \\
& \hspace*{3mm} =
  \bigg\| \sum_{i=1}^m \sum_{j=1}^{n_2}
          \alpha_{ij} (x_{ij} - p) -
          \sum_{i=1}^m \sum_{j=1}^{n_1}
          \alpha_{ij} (x_{ij} - p)
  \bigg\|_a < \varepsilon
\end{split}
\end{eqnarray*}
for all integers $n_1, n_2 > N_\varepsilon$, which implies that
$\big\{ \sum\limits_{i=1}^m \sum\limits_{j=1}^n
        \alpha_{ij} (T_{-q} \circ f \circ T_p)(x_{ij} - p)
 \big\}_n$ is also a Cauchy sequence in $M_a$.
As we proved in \cite[Theorem 2.1]{jungkim} or by Remark
\ref{rem:2.1} $(i)$, we observe that
$(M_a, \langle \cdot, \cdot \rangle_a)$ is a real Hilbert
space.
Thus, $M_a$ is not only complete, but also a Hausdorff space,
so the Cauchy sequence
$\big\{ \sum\limits_{i=1}^m \sum\limits_{j=1}^n
        \alpha_{ij} (T_{-q} \circ f \circ T_p)(x_{ij} - p)
 \big\}_n$ converges in $M_a$, \textit{i.e.}, by Definition
\ref{def:ext}, we have
\begin{eqnarray*}
\begin{split}
(T_{-q} \circ F \circ T_p)
\Bigg( \sum_{i=1}^m \sum_{j=1}^\infty \alpha_{ij} (x_{ij} - p)
\Bigg)
& =   \sum_{i=1}^m \sum_{j=1}^\infty \alpha_{ij}
      (T_{-q} \circ f \circ T_p)(x_{ij} - p) \\
& =   \lim_{n \to \infty}
      \sum_{i=1}^m \sum_{j=1}^n
      \alpha_{ij} (T_{-q} \circ f \circ T_p)(x_{ij} - p) \\
& \in M_a,
\end{split}
\end{eqnarray*}
which implies
\begin{eqnarray*}
F \Bigg( p + \sum_{i=1}^m \sum_{j=1}^\infty
         \alpha_{ij} (x_{ij} - p)
  \Bigg) \in M_a + q = M_a
\end{eqnarray*}
for all $x_{ij} \in E_1$ and $\alpha_{ij} \in \mathbb{R}$ with
$\sum\limits_{i=1}^m \sum\limits_{j=1}^\infty
 \alpha_{ij} (x_{ij} - p) \in M_a$, \textit{i.e.}, the image
of each element of ${\rm GS}(E_1,p)$ under $F$ belongs to
$M_a$.

We now assume that
$\sum\limits_{i=1}^{m_1} \sum\limits_{j=1}^\infty
 \alpha_{ij} (x_{ij} - p)
 = \sum\limits_{i=1}^{m_2} \sum\limits_{j=1}^\infty
   \beta_{ij} (y_{ij} - p) \in M_a$ for some
$m_1, m_2 \in \mathbb{N}$, $x_{ij}, y_{ij} \in E_1$, and
for some $\alpha_{ij}, \beta_{ij} \in \mathbb{R}$.
It then follows from Definition \ref{def:ext} and Lemma
\ref{lem:new2} that
\begin{eqnarray*}
\begin{split}
& \bigg\|
  \big( T_{-q} \circ F \circ T_p \big)
  \bigg( \sum_{i=1}^{m_1} \sum_{j=1}^\infty
         \alpha_{ij} (x_{ij} - p)
  \bigg) -
  \big( T_{-q} \circ F \circ T_p \big)
  \bigg( \sum_{i=1}^{m_2} \sum_{j=1}^\infty
         \beta_{ij} (y_{ij} - p)
  \bigg)
  \bigg\|_a^2 \\
& \hspace*{3mm}
  = \bigg\|
    \sum_{i=1}^{m_1} \sum_{j=1}^\infty \alpha_{ij}
    \big( T_{-q} \circ f \circ T_p \big) (x_{ij} - p) -
    \sum_{i=1}^{m_2} \sum_{j=1}^\infty \beta_{ij}
    \big( T_{-q} \circ f \circ T_p \big) (y_{ij} - p)
    \bigg\|_a^2 \\
& \hspace*{3mm}
  = \bigg\langle
    \sum_{i=1}^{m_1} \sum_{j=1}^\infty \alpha_{ij}
    \big( T_{-q} \circ f \circ T_p \big) (x_{ij} - p) -
    \sum_{i=1}^{m_2} \sum_{j=1}^\infty \beta_{ij}
    \big( T_{-q} \circ f \circ T_p \big) (y_{ij} - p), \\
& \hspace*{11mm}
    \sum_{k=1}^{m_1} \sum_{\ell=1}^\infty \alpha_{k\ell}
    \big( T_{-q} \circ f \circ T_p \big) (x_{k\ell} - p) -
    \sum_{k=1}^{m_2} \sum_{\ell=1}^\infty \beta_{k\ell}
    \big( T_{-q} \circ f \circ T_p \big) (y_{k\ell} - p)
    \bigg\rangle_a \\
& \hspace*{3mm}
  = \cdots \\
& \hspace*{3mm}
  = \bigg\langle
    \sum_{i=1}^{m_1} \sum_{j=1}^\infty \alpha_{ij}
    (x_{ij} - p) -
    \sum_{i=1}^{m_2} \sum_{j=1}^\infty \beta_{ij}
    (y_{ij} - p),\, \\
& \hspace*{11mm}
    \sum_{k=1}^{m_1} \sum_{\ell=1}^\infty \alpha_{k\ell}
    (x_{k\ell} - p) -
    \sum_{k=1}^{m_2} \sum_{\ell=1}^\infty \beta_{k\ell}
    (y_{k\ell} - p)
    \bigg\rangle_a \\
& \hspace*{3mm}
  = \bigg\|
    \sum_{i=1}^{m_1} \sum_{j=1}^\infty \alpha_{ij} (x_{ij} - p) -
    \sum_{i=1}^{m_2} \sum_{j=1}^\infty \beta_{ij} (y_{ij} - p)
    \bigg\|_a^2 \\
& \hspace*{3mm}
  = 0,
\end{split}
\end{eqnarray*}
which implies that
\begin{eqnarray*}
\big( T_{-q} \circ F \circ T_p \big)
\bigg( \sum_{i=1}^{m_1} \sum_{j=1}^\infty
       \alpha_{ij} (x_{ij} - p)
\bigg)
= \big( T_{-q} \circ F \circ T_p \big)
  \bigg( \sum_{i=1}^{m_2} \sum_{j=1}^\infty
         \beta_{ij} (y_{ij} - p)
  \bigg)
\end{eqnarray*}
for all $m_1, m_2 \in \mathbb{N}$, $x_{ij}, y_{ij} \in E_1$,
and for all $\alpha_{ij}, \beta_{ij} \in \mathbb{R}$ satisfying
$\sum\limits_{i=1}^{m_1} \sum\limits_{j=1}^\infty
 \alpha_{ij} (x_{ij} - p)
 = \sum\limits_{i=1}^{m_2} \sum\limits_{j=1}^\infty
   \beta_{ij} (y_{ij} - p) \in M_a$.
\hfill$\Box$
\vspace*{5mm}

We prove in the following theorem that the domain of
$d_a$-isometry $f : E_1 \to E_2$ can be extended to the
first-order generalized linear span ${\rm GS}(E_1,p)$ whenever
$E_1$ is a nonempty bounded subset of $M_a$, whether degenerate
or non-degenerate.
Therefore, Theorem \ref{thm:new3} is a generalization of
\cite[Theorem 2.2]{jung1} for $M_a$.

In the proof, we use the fact that ${\rm GS}(E_1,p) - p$ is a
real vector space.
This fact is self-evident, as briefly mentioned earlier.


\begin{thm} \label{thm:new3}
Assume that $E_1$ and $E_2$ are bounded subsets of $M_a$ that
are $d_a$-isometric to each other via a surjective $d_a$-isometry
$f : E_1 \to E_2$.
Assume that $p$ is an element of $E_1$ and $q$ is an element
of $E_2$ with $q = f(p)$.
The function $F : {\rm GS}(E_1,p) \to M_a$ defined in
Definition \ref{def:ext} is a $d_a$-isometry and the function
$T_{-q} \circ F \circ T_p : {\rm GS}(E_1,p) - p \to M_a$ is a
linear $d_a$-isometry.
In particular, $F$ is an extension of $f$.
\end{thm}

\noindent
\emph{Proof}.
$(a)$
Let $u$ and $v$ be arbitrary elements of the first-order
generalized linear span ${\rm GS}(E_1,p)$ of $E_1$ with
respect to $p$.
Then
\begin{eqnarray} \label{eq:2018-8-23-1}
u - p = \sum_{i=1}^m \sum_{j=1}^\infty
        \alpha_{ij} (x_{ij} - p) \in M_a
~~~\mbox{and}~~~
v - p = \sum_{i=1}^n \sum_{j=1}^\infty
        \beta_{ij} (y_{ij} - p) \in M_a
\end{eqnarray}
for some $m, n \in \mathbb{N}$, some
$x_{ij}, y_{ij} \in E_1$, and for some
$\alpha_{ij}, \beta_{ij} \in \mathbb{R}$.
Then, according to Definition \ref{def:ext}, we have
\begin{eqnarray} \label{eq:2018-8-23-2}
\begin{split}
(T_{-q} \circ F \circ T_p)(u-p)
& = \sum_{i=1}^m \sum_{j=1}^\infty
    \alpha_{ij} (T_{-q} \circ f \circ T_p)(x_{ij} - p), \\
(T_{-q} \circ F \circ T_p)(v-p)
& = \sum_{i=1}^n \sum_{j=1}^\infty
    \beta_{ij} (T_{-q} \circ f \circ T_p)(y_{ij} - p).
\end{split}
\end{eqnarray}

$(b)$
By Lemma \ref{lem:new2}, (\ref{eq:2018-8-23-1}), and
(\ref{eq:2018-8-23-2}), we get
\begin{eqnarray} \label{eq:2018-9-28}
\begin{split}
& \big\langle ( T_{-q} \circ F \circ T_p )(u-p),\,
              ( T_{-q} \circ F \circ T_p )(v-p)
  \big\rangle_a \\
& \hspace*{3mm}
  = \bigg\langle \,\sum_{i=1}^m \sum_{j=1}^\infty \alpha_{ij}
                 (T_{-q} \circ f \circ T_p)(x_{ij} - p),\,
                 \sum_{k=1}^n \sum_{\ell=1}^\infty \beta_{k\ell}
                 (T_{-q} \circ f \circ T_p)(y_{k\ell} - p)
    \bigg\rangle_a \\
& \hspace*{3mm}
  = \sum_{i=1}^m \sum_{k=1}^n
    \sum_{j=1}^\infty \alpha_{ij}
    \sum_{\ell=1}^\infty \beta_{k\ell}
    \big\langle (T_{-q} \circ f \circ T_p)(x_{ij} - p),\,
                (T_{-q} \circ f \circ T_p)(y_{k\ell} - p)
    \big\rangle_a \\
& \hspace*{3mm}
  = \sum_{i=1}^m \sum_{k=1}^n
    \sum_{j=1}^\infty \alpha_{ij}
    \sum_{\ell=1}^\infty \beta_{k\ell}
    \langle x_{ij} - p, y_{k\ell} - p \rangle_a \\
& \hspace*{3mm}
  = \bigg\langle \,\sum_{i=1}^m \sum_{j=1}^\infty \alpha_{ij}
                 (x_{ij} - p),\,
                 \sum_{k=1}^n \sum_{\ell=1}^\infty \beta_{k\ell}
                 (y_{k\ell} - p)
    \bigg\rangle_a \\
& \hspace*{3mm}
  = \langle u - p, v - p \rangle_a
\end{split}
\end{eqnarray}
for all $u, v \in {\rm GS}(E_1,p)$.
That is, $T_{-q} \circ F \circ T_p$ preserves the inner product.
Indeed, equality (\ref{eq:2018-9-28}) is an extended version
of Lemma \ref{lem:new2}.

$(c)$
By using equality (\ref{eq:2018-9-28}), we further obtain
\begin{eqnarray*}
\begin{split}
d_a \big( F(u), F(v) \big)^2
& = \| F(u) - F(v) \|_a^2 \\
& = \big\| ( T_{-q} \circ F \circ T_p )(u-p) -
           ( T_{-q} \circ F \circ T_p )(v-p)
    \big\|_a^2 \\
& = \big\langle ( T_{-q} \circ F \circ T_p )(u-p) -
                ( T_{-q} \circ F \circ T_p )(v-p), \\
&   ~~~~~~      ( T_{-q} \circ F \circ T_p )(u-p) -
                ( T_{-q} \circ F \circ T_p )(v-p)
    \big\rangle_a \\
& = \langle u-p, u-p \rangle_a - \langle u-p, v-p \rangle_a
    - \langle v-p, u-p \rangle_a + \langle v-p, v-p \rangle_a \\
& = \big\langle (u-p) - (v-p), (u-p) - (v-p) \big\rangle_a \\
& = \| (u-p) - (v-p) \|_a^2 \\
& = \| u - v \|_a^2 \\
& = d_a(u, v)^2
\end{split}
\end{eqnarray*}
for all $u, v \in {\rm GS}(E_1,p)$, \textit{i.e.}, $F$ is a $d_a$-isometry.

$(d)$
Now, let $u$ and $v$ be arbitrary elements of ${\rm GS}(E_1,p)$.
Then, it holds that $u - p \in {\rm GS}(E_1,p) - p$,
$v - p \in {\rm GS}(E_1,p) - p$, and
$\alpha (u - p) + \beta (v - p) \in {\rm GS}(E_1,p) - p$ for
any $\alpha, \beta \in \mathbb{R}$, because
${\rm GS}(E_1,p) - p$ is a real vector space.

We get
\begin{eqnarray*}
\begin{split}
& \big\| (T_{-q} \circ F \circ T_p)
         \left( \alpha (u - p) + \beta (v - p) \right) \\
& \hspace*{3mm}
         -\alpha (T_{-q} \circ F \circ T_p) (u - p)
         -\beta (T_{-q} \circ F \circ T_p) (v - p)
  \big\|_a^2 \\
& \hspace*{3mm}
  = \Big\langle
    (T_{-q} \circ F \circ T_p)
    \left( \alpha (u - p) + \beta (v - p) \right) \\
& \hspace*{11mm}
    -\alpha (T_{-q} \circ F \circ T_p) (u - p)
    -\beta (T_{-q} \circ F \circ T_p) (v - p), \\
& \hspace*{11mm}
    (T_{-q} \circ F \circ T_p)
    \left( \alpha (u - p) + \beta (v - p) \right) \\
& \hspace*{11mm}
    -\alpha (T_{-q} \circ F \circ T_p) (u - p)
    -\beta (T_{-q} \circ F \circ T_p) (v - p)
    \Big\rangle_a. 
\end{split}
\end{eqnarray*}
Since $\alpha (u - p) + \beta (v - p) = w - p$ for some
$w \in {\rm GS}(E_1,p)$, we further use (\ref{eq:2018-9-28})
to obtain
\begin{eqnarray*}
\begin{split}
& \big\| (T_{-q} \circ F \circ T_p)
         \left( \alpha (u - p) + \beta (v - p) \right) \\
& \hspace*{3mm}
         -\alpha (T_{-q} \circ F \circ T_p) (u - p)
         -\beta (T_{-q} \circ F \circ T_p) (v - p)
  \big\|_a^2 \\
& \hspace*{3mm}
  = \langle w - p,\, w - p \rangle_a -
    \alpha \langle w - p,\, u - p \rangle_a -
    \beta \langle w - p,\, v - p \rangle_a \\
& \hspace*{8mm}
    -\alpha \langle u - p,\, w - p \rangle_a
    +\alpha^2 \langle u - p,\, u - p \rangle_a
    +\alpha \beta \langle u - p,\, v - p \rangle_a \\
& \hspace*{8mm}
    -\beta \langle v - p,\, w - p \rangle_a
    +\alpha \beta \langle v - p,\, u - p \rangle_a
    +\beta^2 \langle v - p,\, v - p \rangle_a \\
& \hspace*{3mm}
  = 0,
\end{split}
\end{eqnarray*}
which implies that the function
$T_{-q} \circ F \circ T_p : {\rm GS}(E_1,p) - p \to M_a$ is
linear.

$(e)$
Finally, we set $\alpha_{11} = 1$, $\alpha_{ij} = 0$ for any
$(i,j) \neq (1,1)$, and $x_{11} = x$ in (\ref{eq:2018-8-23-1})
and (\ref{eq:2018-8-23-2}) to see
\begin{eqnarray*}
(T_{-q} \circ F \circ T_p)(x-p)
= (T_{-q} \circ f \circ T_p)(x-p)
\end{eqnarray*}
for every $x \in E_1$, which implies that $F(x) = f(x)$ for
every $x \in E_1$, \textit{i.e.}, $F$ is an extension of $f$.
\hfill$\Box$


\section{Second-order generalized linear span}
\label{sec:5}

For any element $x$ of $M_a$ and $r > 0$, we denote by
$B_r(x)$ the open ball defined by
$B_r(x) = \{ y \in M_a : \| y - x \|_a < r \}$.

Definitions \ref{def:gls} and \ref{def:ext} will be generalized
to the cases of $n > 1$ in the following definition.
We introduce the concept of $n$th-order generalized linear
spans ${\rm GS}^n(E_1,p)$, which generalizes the concept of
first-order generalized linear span ${\rm GS}(E,p)$.
Moreover, we define the $d_a$-isometries $F_n$ which extend
the domain of $d_a$-isometry $f$ to ${\rm GS}^n(E_1,p)$.

It is surprising, however, that this process of generalization
does not go far.
Indeed, we will find in Remark \ref{rem:new3.1} and Theorem
\ref{thm:setyol} that ${\rm GS}^2(E_1,p)$ and $F_2$ are their
limits.


\begin{defn} \label{def:6.1}
Let $E_1$ be a nonempty bounded subset of $M_a$ that is
$d_a$-isometric to a subset $E_2$ of $M_a$ via a surjective
$d_a$-isometry $f : E_1 \to E_2$.
Let $p$ be an element of $E_1$ and $q$ an element of $E_2$
with $q = f(p)$.
Assume that $r$ is a positive real number satisfying
$E_1 \subset B_r(p)$.
\begin{itemize}
\item[$(i)$]  We define ${\rm GS}^0(E_1,p) = E_1$ and
              ${\rm GS}^1(E_1,p) = {\rm GS}(E_1,p)$.
              In general, we define the
              \emph{$n$th-order generalized linear span} of
              $E_1$ with respect to $p$ as
              ${\rm GS}^n(E_1,p) = {\rm GS}
               ({\rm GS}^{n-1}(E_1,p) \cap B_r(p), p)$ for all
              $n \in \mathbb{N}$.
\item[$(ii)$] We define $F_0 = f$ and $F_1 = F$, where $F$ is
              defined in Definition \ref{def:ext}.
              Moreover, for any $n \in \mathbb{N}$, we define
              the function $F_n : {\rm GS}^n(E_1,p) \to M_a$
              by
              \begin{eqnarray*}
              (T_{-q} \circ F_n \circ T_p)
              \bigg( \sum_{i=1}^m \sum_{j=1}^\infty
                     \alpha_{ij} (x_{ij} - p) \bigg)
              = \sum_{i=1}^m \sum_{j=1}^\infty \alpha_{ij}
                (T_{-q} \circ F_{n-1} \circ T_p) (x_{ij} - p)
              \end{eqnarray*}
              for all $m \in \mathbb{N}$,
              $x_{ij} \in {\rm GS}^{n-1}(E_1,p) \cap B_r(p)$,
              and $\alpha_{ij} \in \mathbb{R}$ with
              $\sum\limits_{i=1}^m \sum\limits_{j=1}^\infty
               \alpha_{ij} (x_{ij} - p) \in M_a$.
\end{itemize}
\end{defn}


\begin{rem} \label{rem:31}
Let $E$ be a nonempty bounded subset of $M_a$.
If $s$ and $t$ are positive real numbers that satisfy
$E \subset B_s(p) \cap B_t(p)$, then
\begin{eqnarray*}
{\rm GS} \big( {\rm GS}(E,p) \cap B_s(p), p \big)
= {\rm GS} \big( {\rm GS}(E,p) \cap B_t(p), p \big).
\end{eqnarray*}
\end{rem}

\noindent
\emph{Proof}.
Assume that $0 < s < t$.
Then, there exists a real number $c > 1$ with $s > \frac{t}{c}$
and it is obvious that $B_{t/c}(p) \subset B_s(p)$.
Assume that $x$ is an arbitrary element of
${\rm GS}({\rm GS}(E,p) \cap B_t(p), p)$.
Then there exist some $m \in \mathbb{N}$, some
$u_{ij} \in {\rm GS}(E,p) \cap B_t(p)$ and some
$\alpha_{ij} \in \mathbb{R}$ such that
$x = p + \sum\limits_{i=1}^m \sum\limits_{j=1}^\infty
 \alpha_{ij} (u_{ij} - p) \in M_a$.
We notice that
\begin{eqnarray*}
\big( {\rm GS}(E,p) - p \big) \cap \big( B_t(p) - p \big)
= \big\{ u - p \in M_a : u \in {\rm GS}(E,p) \cap B_t(p) \big\}.
\end{eqnarray*}

Since ${\rm GS}(E,p) - p$ is a real vector space,
$\frac{t}{c} < s$, and since
$u_{ij} - p \in ({\rm GS}(E,p) - p) \cap (B_t(p) - p)$ for
any $i$ and $j$, we have
\begin{eqnarray*}
\frac{1}{c} (u_{ij} - p)
\in ({\rm GS}(E,p) - p) \cap (B_s(p) - p).
\end{eqnarray*}
Hence, we can choose a $v_{ij} \in {\rm GS}(E,p) \cap B_s(p)$
such that $\frac{1}{c} (u_{ij} - p) = v_{ij} - p$.
Thus, we get
\begin{eqnarray*}
x =   p + \sum_{i=1}^m \sum_{j=1}^\infty
      \alpha_{ij} (u_{ij} - p)
  =   p + \sum_{i=1}^m \sum_{j=1}^\infty
      c \alpha_{ij} (v_{ij} - p)
  \in {\rm GS}({\rm GS}(E,p) \cap B_s(p), p),
\end{eqnarray*}
which implies that
${\rm GS}({\rm GS}(E,p) \cap B_t(p), p)
 \subset {\rm GS}({\rm GS}(E,p) \cap B_s(p), p)$.

The reverse inclusion is obvious, since
$B_s(p) \subset B_t(p)$.
\hfill$\Box$
\vspace*{5mm}

We generalize Lemma \ref{lem:new2} and formula
(\ref{eq:2018-9-28}) in the following lemma.
Indeed, we prove that the function
$T_{-q} \circ F_n \circ T_p : {\rm GS}^n(E_1,p) - p \to M_a$
preserves the inner product.
This property is important in proving the following theorems
as a necessary condition for $F_n$ to be a $d_a$-isometry.


\begin{lem} \label{lem:6.1}
Let $E_1$ be a bounded subset of $M_a$ that is $d_a$-isometric to
a subset $E_2$ of $M_a$ via a surjective $d_a$-isometry
$f : E_1 \to E_2$.
Assume that $p$ and $q$ are elements of $E_1$ and $E_2$, which
satisfy $q = f(p)$.
If $n \in \mathbb{N}$, then
\begin{eqnarray*}
\big\langle (T_{-q} \circ F_n \circ T_p) (u - p),
            (T_{-q} \circ F_n \circ T_p) (v - p)
\big\rangle_a =
\langle u - p, v - p \rangle_a
\end{eqnarray*}
for all $u, v \in {\rm GS}^n(E_1,p)$.
\end{lem}

\noindent
\emph{Proof}.
Our assertion for $n = 1$ was already proved in
(\ref{eq:2018-9-28}).
Considering Remark \ref{rem:31}, assume that $r$ is a positive
real number satisfying $E_1 \subset B_r(p)$.
Now we assume that the assertion is true for some
$n \in \mathbb{N}$.
Let $u, v$ be arbitrary elements of ${\rm GS}^{n+1}(E_1,p)$.
Then there exist some $m_1, m_2 \in \mathbb{N}$, some
$x_{ij}, y_{k\ell} \in {\rm GS}^n(E_1,p) \cap B_r(p)$ and some
$\alpha_{ij}, \beta_{k\ell} \in \mathbb{R}$ such that
\begin{eqnarray*}
u - p = \sum_{i=1}^{m_1} \sum_{j=1}^\infty
        \alpha_{ij} (x_{ij} - p) \in M_a ~~~\mbox{and}~~~
v - p = \sum_{k=1}^{m_2} \sum_{\ell=1}^\infty
        \beta_{k\ell} (y_{k\ell} - p) \in M_a.
\end{eqnarray*}

Using Definition \ref{def:6.1} $(ii)$ and our assumption, we
get
\begin{eqnarray*}
\begin{split}
& \big\langle ( T_{-q} \circ F_{n+1} \circ T_p )(u - p),\,
              ( T_{-q} \circ F_{n+1} \circ T_p )(v - p)
  \big\rangle_a \\
& \hspace*{3mm}
  = \bigg\langle \,\sum_{i=1}^{m_1} \sum_{j=1}^\infty
                 \alpha_{ij}
                 (T_{-q} \circ F_n \circ T_p)(x_{ij} - p),\,
                 \sum_{k=1}^{m_2} \sum_{\ell=1}^\infty
                 \beta_{k\ell}
                 (T_{-q} \circ F_n \circ T_p)(y_{k\ell} - p)
    \bigg\rangle_a \\
& \hspace*{3mm}
  = \sum_{i=1}^{m_1} \sum_{k=1}^{m_2}
    \sum_{j=1}^\infty \alpha_{ij}
    \sum_{\ell=1}^\infty \beta_{k\ell}
    \big\langle (T_{-q} \circ F_n \circ T_p)(x_{ij} - p),\,
                (T_{-q} \circ F_n \circ T_p)(y_{k\ell} - p)
    \big\rangle_a \\
& \hspace*{3mm}
  = \sum_{i=1}^{m_1} \sum_{k=1}^{m_2}
    \sum_{j=1}^\infty \alpha_{ij}
    \sum_{\ell=1}^\infty \beta_{k\ell}
    \langle x_{ij} - p, y_{k\ell} - p \rangle_a \\
& \hspace*{3mm}
  = \bigg\langle \,\sum_{i=1}^{m_1} \sum_{j=1}^\infty
                 \alpha_{ij} (x_{ij} - p),\,
                 \sum_{k=1}^{m_2} \sum_{\ell=1}^\infty
                 \beta_{k\ell} (y_{k\ell} - p)
    \bigg\rangle_a \\
& \hspace*{3mm}
  = \langle u - p, v - p \rangle_a
\end{split}
\end{eqnarray*}
for all $u, v \in {\rm GS}^{n+1}(E_1,p)$.
By mathematical induction, we may then conclude that our
assertion is true for all $n \in \mathbb{N}$.
\hfill$\Box$
\vspace*{5mm}

When $n = 1$ and $p = p'$, the first assertion in $(i)$ of the
following lemma is self-evident, so we have used that fact
several times before, omitting the proof.
The assertion $(iv)$ in the following lemma seems to be related
in some way to Remark \ref{rem:31}.


\begin{lem} \label{lem:3.8}
Let $E$ be a bounded subset of $M_a$ and $p, p' \in E$.
Assume that $r$ is a positive real number satisfying
$E \subset B_r(p)$.
\begin{itemize}
\item[$(i)$]   ${\rm GS}^n(E,p) - p'$ is a vector space over
               $\mathbb{R}$ for each $n \in \mathbb{N}$.
\item[$(ii)$]  ${\rm GS}^n(E,p) \subset {\rm GS}^{n+1}(E,p)$
               for each $n \in \mathbb{N}$.
\item[$(iii)$] ${\rm GS}^2(E,p) = \overline{{\rm GS}(E,p)}$,
               where $\overline{{\rm GS}(E,p)}$ is the closure
               of\/ ${\rm GS}(E,p)$ in $M_a$.
\item[$(iv)$]  $\Lambda({\rm GS}^n(E,p)) =
                \Lambda({\rm GS}^n(E,p) \cap B_r(p))$ for all
               $n \in \mathbb{N}$.
\end{itemize}
\end{lem}

\noindent
\emph{Proof}.
$(i)$ By using Definitions \ref{def:gls} and \ref{def:6.1}, we
prove that ${\rm GS}(E,p) - p'$ is a real vector space.
(We can prove similarly for the case of $n > 1$.)
Given $x, y \in {\rm GS}(E,p) - p'$, we may choose some
$m_1, m_2 \in \mathbb{N}$, some $u_{ij}, v_{ij} \in E$, and
some $\alpha_{ij}, \beta_{ij} \in \mathbb{R}$ such that
$x = (p - p') + \sum\limits_{i=1}^{m_1} \sum\limits_{j=1}^\infty
     \alpha_{ij} (u_{ij} - p) \in M_a$ and
$y = (p - p') + \sum\limits_{i=1}^{m_2} \sum\limits_{j=1}^\infty
     \beta_{ij} (v_{ij} - p) \in M_a$.
Since $M_a$ is a real vector space,
$\alpha \sum\limits_{i=1}^{m_1} \sum\limits_{j=1}^\infty
 \alpha_{ij} (u_{ij} - p) +
 \beta \sum\limits_{i=1}^{m_2} \sum\limits_{j=1}^\infty
 \beta_{ij} (v_{ij} - p) \in M_a$ for all
$\alpha, \beta \in \mathbb{R}$.

Moreover, we see that
\begin{eqnarray*}
\begin{split}
\alpha x + \beta y
& =   \left(
      p + (1 - \alpha - \beta)(p' - p) +
      \sum_{i=1}^{m_1} \sum_{j=1}^\infty
      \alpha \alpha_{ij} (u_{ij} - p) +
      \sum_{i=1}^{m_2} \sum_{j=1}^\infty
      \beta \beta_{ij} (v_{ij} - p)
      \right) - p' \\
& \in {\rm GS}(E,p) - p'
\end{split}
\end{eqnarray*}
for all $\alpha, \beta \in \mathbb{R}$.
Hence, ${\rm GS}(E,p) - p'$ is a real vector space as a
subspace of real vector space $M_a$.

$(ii)$
Let $r$ be a positive real number with $E \subset B_r(p)$.
If $x \in {\rm GS}^n(E,p)$ for some $n \in \mathbb{N}$, then
$x - p \in {\rm GS}^n(E,p) - p$.
Since ${\rm GS}^n(E,p) - p$ is a real vector space by $(i)$
and $B_r(p) - p = B_r(0)$, we can choose a (positive or
negative but sufficiently small) real number $\mu \neq 0$ such
that $\mu (x - p) \in ({\rm GS}^n(E,p) - p) \cap (B_r(p) - p)$.
We notice that
\begin{eqnarray} \label{eq:2021-11-2}
\big( {\rm GS}^n(E,p) - p \big) \cap \big( B_r(p) - p \big)
= \big\{ v - p \in M_a : v \in {\rm GS}^n(E,p) \cap B_r(p)
  \big\}.
\end{eqnarray}
Thus, we see that $\mu (x - p) = v - p$ for some
$v \in {\rm GS}^n(E,p) \cap B_r(p)$.
Since $x = p + \frac{1}{\mu} (v - p)$, it holds that
$x \in {\rm GS}^{n+1}(E,p)$.
Therefore, we conclude that
${\rm GS}^n(E,p) \subset {\rm GS}^{n+1}(E,p)$ for every
$n \in \mathbb{N}$.

$(iii)$ Let $x$ be an arbitrary element of
$\overline{{\rm GS}(E,p)}$.
Then there exists some sequence $\{ x_n \}$ in ${\rm GS}(E,p)$
that converges to $x$, where $x_n \neq x$ for all
$n \in \mathbb{N}$.
We now set $y_1 = x_1$ and $y_i = x_i - x_{i-1}$ for each
integer $i \geq 2$.
Then we have
\begin{eqnarray*}
x_n = \sum_{i=1}^n y_i,
\end{eqnarray*}
where $y_i = (x_i - p) - (x_{i-1} - p) \in {\rm GS}(E,p) - p$
for $i \geq 2$.
Since ${\rm GS}(E,p) - p$ is a real vector space and
$B_r(p) - p = B_r(0)$, we can select a real number
$\mu_i \neq 0$ such that
\begin{eqnarray*}
\mu_i y_i \in {\rm GS}(E,p) - p ~~~\mbox{and}~~~
\mu_i y_i \in B_r(p) - p
\end{eqnarray*}
for every integer $i \geq 2$.
Thus, it follows from (\ref{eq:2021-11-2}) that
\begin{eqnarray*}
x_n = \sum_{i=1}^n y_i
    = y_1 + \sum_{i=2}^n \frac{1}{\mu_i} (\mu_i y_i)
    = x_1 + \sum_{i=2}^n \frac{1}{\mu_i} (v_i - p),
\end{eqnarray*}
where $v_i \in {\rm GS}(E,p) \cap B_r(p)$ for $i \geq 2$.
Since the sequence $\{ x_n \}$ is assumed to converge to $x$,
the sequence
$\big\{ x_1 + \sum\limits_{i=2}^n \frac{1}{\mu_i} (v_i - p)
 \big\}_n$ converges to $x$.
Hence, we have
\begin{eqnarray} \label{eq:2021-10-16-1}
x_1 + \sum_{i=2}^\infty \frac{1}{\mu_i} (v_i - p)
=   \lim_{n \to \infty} x_n
=   x
\in M_a.
\end{eqnarray}
(Since $M_a$ is a Huasdorff space, $x$ is the unique limit
point of the sequence $\{ x_n \}$.)

Furthermore, there exists a real number $\mu_1 \neq 0$ that
satisfies $\mu_1 (x_1 - p) \in {\rm GS}(E,p) - p$ and
$\mu_1 (x_1 - p) \in B_r(p) - p$, \textit{i.e.},
$\mu_1 (x_1 - p) \in ({\rm GS}(E,p) - p) \cap (B_r(p) - p)$.
Thus, there exists a $v_1 \in {\rm GS}(E,p) \cap B_r(p)$ such
that $\mu_1 (x_1 - p) = v_1 - p$ or
$x_1 - p = \frac{1}{\mu_1} (v_1 - p)$.
Therefore,
\begin{eqnarray} \label{eq:2021-10-16-2}
x =   p + (x_1 - p) + \sum_{i=2}^\infty
      \frac{1}{\mu_i} (v_i - p)
  =   p + \sum_{i=1}^\infty \frac{1}{\mu_i} (v_i - p),
\end{eqnarray}
where $v_i \in {\rm GS}(E,p) \cap B_r(p)$ for each
$i \in \mathbb{N}$.
It follows from (\ref{eq:2021-10-16-1}) that
$\sum\limits_{i=1}^\infty \frac{1}{\mu_i} (v_i - p) \in M_a$.
Thus, by (\ref{eq:2021-10-16-2}), we see that
$x \in {\rm GS}^2(E,p)$, which implies that
$\overline{{\rm GS}(E,p)} \subset {\rm GS}^2(E,p)$.

On the other hand, let $y \in {\rm GS}^2(E,p)$.
Then there are some $m \in \mathbb{N}$, some
$v_{ij} \in {\rm GS}(E,p) \cap B_r(p)$, and some
$\alpha_{ij} \in \mathbb{R}$ such that
$y = p + \sum\limits_{i=1}^m \sum\limits_{j=1}^\infty
 \alpha_{ij} (v_{ij} - p) \in M_a$.
Let us define
$y_n = p + \sum\limits_{i=1}^m \sum\limits_{j=1}^n
 \alpha_{ij} (v_{ij} - p)$ for every $n \in \mathbb{N}$.
Since $v_{ij} - p \in {\rm GS}(E,p) - p$ for all $i$ and $j$
and ${\rm GS}(E,p) - p$ is a real vector space, we know that
$y_n - p = \sum\limits_{i=1}^m \sum\limits_{j=1}^n
 \alpha_{ij} (v_{ij} - p) \in {\rm GS}(E,p) - p$ and hence,
$y_n \in {\rm GS}(E,p)$ for all $n \in \mathbb{N}$.
Since ${\rm GS}(E,p)$ is a Hausdorff space, $y$ is the unique
limit point of the sequence $\{ y_n \}$.
Thus, we see that
\begin{eqnarray*}
y =   p + \sum_{i=1}^m \sum_{j=1}^\infty
      \alpha_{ij} (v_{ij} - p)
  =   \lim_{n \to \infty} y_n
  \in \overline{{\rm GS}(E,p)},
\end{eqnarray*}
which implies that
${\rm GS}^2(E,p) \subset \overline{{\rm GS}(E,p)}$.

$(iv)$ Let $i \in \Lambda({\rm GS}^n(E,p))$.
In view of Definition \ref{def:Lambda}, there exist
$x \in {\rm GS}^n(E,p)$ and $\alpha \neq 0$ with
$x + \alpha e_i \in {\rm GS}^n(E,p)$.
Further,
$x = p + \sum\limits_{i=1}^m \sum\limits_{j=1}^\infty
 \alpha_{ij} (u_{ij} - p)$ for some $m \in \mathbb{N}$, some
$u_{ij} \in {\rm GS}^{n-1}(E,p) \cap B_r(p)$, and for some
$\alpha_{ij} \in \mathbb{R}$.
Since
$\sum\limits_{i=1}^m \sum\limits_{j=1}^\infty
 \alpha_{ij} (u_{ij} - p) + \alpha e_i
 =   x - p + \alpha e_i
 \in {\rm GS}^n(E,p) - p$, it holds that
$\mu \big( \sum\limits_{i=1}^m \sum\limits_{j=1}^\infty
           \alpha_{ij} (u_{ij} - p) + \alpha e_i \big)
 \in ({\rm GS}^n(E,p) - p) \cap (B_r(p) - p)$
for any sufficiently small $\mu \neq 0$, or equivalently, it
follows from (\ref{eq:2021-11-2}) that
\begin{eqnarray} \label{eq:2021-4-28}
\left( p + \sum\limits_{i=1}^m \sum\limits_{j=1}^\infty
       \mu \alpha_{ij} (u_{ij} - p)
\right) + \mu \alpha e_i
\in {\rm GS}^n(E,p) \cap B_r(p).
\end{eqnarray}

On the other hand, since
$\sum\limits_{i=1}^m \sum\limits_{j=1}^\infty
 \alpha_{ij} (u_{ij} - p) = x - p
 \in {\rm GS}^n(E,p) - p$, it holds that
$\sum\limits_{i=1}^m \sum\limits_{j=1}^\infty
 \mu \alpha_{ij} (u_{ij} - p) \in ({\rm GS}^n(E,p) - p)
 \cap (B_r(p) - p)$ for any sufficiently small $\mu \neq 0$.
Hence, it follows from (\ref{eq:2021-11-2}) that
$p + \sum\limits_{i=1}^m \sum\limits_{j=1}^\infty
 \mu \alpha_{ij} (u_{ij} - p)  \in {\rm GS}^n(E,p)
 \cap B_r(p)$ for any sufficiently small $\mu \neq 0$.
Thus, by Definition \ref{def:Lambda} and (\ref{eq:2021-4-28}),
it holds that $i \in \Lambda({\rm GS}^n(E,p) \cap B_r(p))$,
which implies that
$\Lambda({\rm GS}^n(E,p))
 \subset \Lambda({\rm GS}^n(E,p) \cap B_r(p))$.
Obviously, the inverse inclusion is true.
\hfill$\Box$
\vspace*{5mm}

As we mentioned earlier, we will see that the second-order
generalized linear span is the last step in this kind of domain
extension.


\begin{rem} \label{rem:new3.1}
If $E$ is a bounded subset of $M_a$ and $p \in E$, then
\begin{eqnarray*}
E \subset {\rm GS}(E,p) \subset \overline{{\rm GS}(E,p)}
= {\rm GS}^2(E,p) = {\rm GS}^n(E,p)
\end{eqnarray*}
for any integer $n \geq 3$.
Indeed, ${\rm GS}^n(E,p) - p$ is a real Hilbert space for
$n \geq 2$.
\end{rem}

\noindent
\emph{Proof}.
$(a)$
Considering Remark \ref{rem:31}, we can choose a real number
$r > 0$ that satisfies $E \subset B_r(p)$.
Assume that $x \in {\rm GS}^3(E,p)$.
Then there exist some $m_0 \in \mathbb{N}$, some
$u_{ij} \in {\rm GS}^2(E,p) \cap B_r(p)$, and some
$\alpha_{ij} \in \mathbb{R}$ such that
$x = p + \sum\limits_{i=1}^{m_0} \sum\limits_{j=1}^\infty
 \alpha_{ij} (u_{ij} - p) \in M_a$.

We define
$x_m = p + \sum\limits_{i=1}^{m_0} \sum\limits_{j=1}^m
 \alpha_{ij} (u_{ij} - p)$ for each $m \in \mathbb{N}$.
Since $u_{ij} \in {\rm GS}^2(E,p)$, there exist some
$m_{ij} \in \mathbb{N}$, some
$v_{ijk\ell} \in {\rm GS}(E,p) \cap B_r(p)$, and some
$\beta_{ijk\ell} \in \mathbb{R}$ such that
$u_{ij} = p + \sum\limits_{k=1}^{m_{ij}}
 \sum\limits_{\ell=1}^\infty
 \beta_{ijk\ell} (v_{ijk\ell} - p) \in M_a$.
Hence, it holds that
\begin{eqnarray*}
x_m =   p + \sum_{i=1}^{m_0} \sum_{j=1}^m \sum_{k=1}^{m_{ij}}
            \sum_{\ell=1}^\infty
            \alpha_{ij} \beta_{ijk\ell}
            \big( v_{ijk\ell} - p \big)
    \in M_a,
\end{eqnarray*}
which implies that $x_m \in {\rm GS}^2(E,p)$ for all
$m \in \mathbb{N}$.
Thus, $\{ x_m \}$ is a sequence in ${\rm GS}^2(E,p)$ that
converges to $x$.
Therefore, $x \in {\rm GS}^2(E,p)$ because ${\rm GS}^2(E,p)$
is closed.
Thus, ${\rm GS}^3(E,p) \subset {\rm GS}^2(E,p)$.
The inverse inclusion is of course true due to Lemma
\ref{lem:3.8} $(ii)$.
We have proved that ${\rm GS}^2(E,p) = {\rm GS}^3(E,p)$.

$(b)$ Assume that
${\rm GS}^2(E,p) = \cdots = {\rm GS}^n(E,p) = {\rm GS}^{n+1}(E,p)$
for some integer $n \geq 2$.

$(c)$ If we replace ${\rm GS}(E,p)$, ${\rm GS}^2(E,p)$, and
${\rm GS}^3(E,p)$ in the previous part $(a)$ with
${\rm GS}^n(E,p)$, ${\rm GS}^{n+1}(E,p)$, and
${\rm GS}^{n+2}(E,p)$, respectively, and if we consider the
fact that ${\rm GS}^{n+1}(E,p) = {\rm GS}^2(E,p)$ is closed in
$M_a$ by Lemma \ref{lem:3.8} $(iii)$ and our assumption $(b)$,
then we arrive at the conclusion that
${\rm GS}^{n+1}(E,p) = {\rm GS}^{n+2}(E,p)$.

$(d)$ With the conclusion of mathematical induction we prove
that ${\rm GS}^n(E,p) = {\rm GS}^2(E,p)$ for every integer
$n \geq 3$.
Moreover, when $n \geq 2$, ${\rm GS}^n(E,p)$ is complete as a
closed subset of a real Hilbert space $M_a$
(ref. Remark \ref{rem:2.1} and \cite[Theorem 63.3]{kasriel}).
Therefore, ${\rm GS}^n(E,p) - p$ is a real Hilbert space for
$n \geq 2$.
\hfill$\Box$
\vspace*{5mm}

The following lemma is an extension of Lemma \ref{lem:2018-9-18}
for the second-order generalized linear span ${\rm GS}^2(E,p)$.
Indeed, we prove that if $i \in \Lambda({\rm GS}^2(E,p))$, then
the second-order generalized linear span of $E$ contains all
the lines through ${\rm GS}(E,p)$ in the direction $e_i$.


\begin{lem} \label{lem:setahob}
Assume that a bounded subset $E$ of $M_a$ contains at least
two elements and $p \in E$.
If $i \in \Lambda({\rm GS}^2(E,p))$ and $p' \in {\rm GS}(E,p)$,
then $p' + \alpha_i e_i \in {\rm GS}^2(E,p)$ for any
$\alpha_i \in \mathbb{R}$.
\end{lem}

\noindent
\emph{Proof}.
Let $r$ be a positive real number with $E \subset B_r(p)$.
Assume that $i \in \Lambda({\rm GS}^2(E,p))$.
Considering Lemma \ref{lem:3.8} $(iv)$ and Remark
\ref{rem:new3.1}, if we substitute
${\rm GS}^2(E,p) \cap B_{r}(p)$ for $E$ in Lemma
\ref{lem:2018-9-18}, then
$p + \alpha_i e_i \in {\rm GS}^3(E,p) = {\rm GS}^2(E,p)$ for
all $\alpha_i \in \mathbb{R}$.
Thus, there are some $m \in \mathbb{N}$, some
$w_{ij} \in {\rm GS}(E,p) \cap B_{r}(p)$, and some
$\beta_{ij} \in \mathbb{R}$ with
$\sum\limits_{i=1}^m \sum\limits_{j=1}^\infty
 \beta_{ij} (w_{ij} - p) \in M_a$ such that
$p + \alpha_i e_i
 = p + \sum\limits_{i=1}^m \sum\limits_{j=1}^\infty
   \beta_{ij} (w_{ij} - p)$, and hence, we have
\begin{eqnarray} \label{eq:2020-1-21-1}
p' + \alpha_i e_i
= p + \alpha_i e_i + (p' - p)
= p + \sum_{i=1}^m \sum_{j=1}^\infty
  \beta_{ij} (w_{ij} - p) + (p' - p).
\end{eqnarray}
Because $p' - p$ belongs to ${\rm GS}(E,p) - p$, which is a
real vector space by Lemma \ref{lem:3.8} $(i)$, and
$B_r(p) - p = B_r(0)$, we can choose some sufficiently small
real number $\mu \neq 0$ such that
\begin{eqnarray} \label{eq:2020-1-21-2}
\mu (p' - p) \in {\rm GS}(E,p) - p ~~~\mbox{and}~~~
\mu (p' - p) \in B_{r}(p) - p.
\end{eqnarray}

Considering (\ref{eq:2021-11-2}), (\ref{eq:2020-1-21-1}) and
(\ref{eq:2020-1-21-2}), if we put $\mu (p' - p) = w - p$ with
a $w \in {\rm GS}(E,p) \cap B_{r}(p)$, then we have
\begin{eqnarray*}
p' + \alpha_i e_i
=   p + \sum_{i=1}^m \sum_{j=1}^\infty
    \beta_{ij} (w_{ij} - p) + \frac{1}{\mu} (w - p)
\in {\rm GS}^2(E,p)
\end{eqnarray*}
for all $\alpha_i \in \mathbb{R}$.
\hfill$\Box$


\section{Basic cylinders and basic intervals}
\label{sec:6}

First, we will define the infinite dimensional intervals,
which were simply defined in \cite{jungkim}, more precisely
divided into non-degenerate basic cylinders, degenerate basic
cylinders, and basic intervals.


\begin{defn} \label{def:basic}
For any positive integer $n$, we define the infinite
dimensional interval by
\begin{eqnarray*}
J = \prod_{i=1}^\infty J_i, ~~~\mbox{where}~~~
J_i = \left\{\begin{array}{ll}
                \mbox{$[ 0, p_{2i} ]$}
                & (\mbox{for}~ i \in \Lambda_1),
                \vspace*{2mm}\\
                \mbox{$[ p_{1i}, p_{2i} ]$}
                & (\mbox{for}~ i \in \Lambda_2),
                \vspace*{2mm}\\
                \mbox{$[ p_{1i}, 1 ]$}
                & (\mbox{for}~ i \in \Lambda_3),
                \vspace*{2mm}\\
                \mbox{$\{ p_{1i} \}$}
                & (\mbox{for}~ i \in \Lambda_4),
                \vspace*{2mm}\\
                \mbox{$[ 0, 1 ]$}
                & (\mbox{otherwise})
             \end{array}
      \right.
\end{eqnarray*}
for some disjoint finite subsets $\Lambda_1$, $\Lambda_2$,
$\Lambda_3$ of $\{ 1, 2, \ldots, n \}$ and
$0 < p_{1i} < p_{2i} < 1$ for
$i \in \Lambda_1 \cup \Lambda_2 \cup \Lambda_3$ and
$0 \leq p_{1i} \leq 1$ for $i \in \Lambda_4$.
If $\Lambda_4 = \emptyset$, then $J$ is called a
\emph{non-degenerate basic cylinder}.
When $\Lambda_4$ is a nonempty finite set, $J$ is called a
\emph{degenerate basic cylinder}.
If $\Lambda_4$ is an infinite set, then $J$ will be called a
\emph{basic interval}.
\end{defn}

We note that in order for an infinite dimensional interval $J$
to become a basic cylinder, $\Lambda_4$ must be a finite set.
Further, we remark that
$\Lambda_4 = \mathbb{N} \!\setminus\! \Lambda(J)$ and
$\Lambda(J) = \mathbb{N} \!\setminus\! \Lambda_4$.
That is, $\mathbb{N}$ is the disjoint union of $\Lambda(J)$
and $\Lambda_4$.


\begin{thm} \label{thm:set9}
Let $J$ be either a translation of a basic cylinder or a
translation of a basic interval and $p \in J$.
Then
\begin{eqnarray*}
{\rm GS}(J,p)
= \Bigg\{ p + \sum_{i \in \Lambda(J)} \alpha_i e_i \in M_a :
          \alpha_i \in \mathbb{R} ~\mbox{for all}~
          i \in \Lambda(J)
  \Bigg\}.
\end{eqnarray*}
\end{thm}

\noindent
\emph{Proof}.
Assume that $x$ is an arbitrary element of ${\rm GS}(J,p)$.
Since $\{ \frac{1}{a_i} e_i \}$ is a complete orthonormal
sequence in $M_a$ and
$\langle x-p, \frac{1}{a_i} e_i \rangle_a = a_i (x_i - p_i)$
for all $i$, it follows from Lemma \ref{lem:doolset}
that
\begin{eqnarray*}
\begin{split}
x - p & = \sum_{i=1}^\infty
          \bigg\langle x-p, \frac{1}{a_i} e_i \bigg\rangle_a
          \frac{1}{a_i} e_i \\
      & = \sum_{i=1}^\infty a_i (x_i - p_i) \frac{1}{a_i} e_i \\
      & = \sum_{i \in \Lambda(J)}
          a_i (x_i - p_i) \frac{1}{a_i} e_i \\
      & = \sum_{i \in \Lambda(J)}
          \bigg\langle x-p, \frac{1}{a_i} e_i \bigg\rangle_a
          \frac{1}{a_i} e_i.
\end{split}
\end{eqnarray*}
Hence, we see that $x \in M_a$ and
\begin{eqnarray*}
\begin{split}
x & =   p + \sum_{i \in \Lambda(J)}
        \left\langle x - p,\, \frac{1}{a_i} e_i \right\rangle_a
        \frac{1}{a_i} e_i ~(\in M_a) \\
  & \in \Bigg\{ p + \sum_{i \in \Lambda(J)}
                \alpha_i e_i \in M_a : \alpha_i \in \mathbb{R}
                ~\mbox{for all}~ i \in \Lambda(J)
        \Bigg\},
\end{split}
\end{eqnarray*}
which implies that
\begin{eqnarray*}
\hspace*{3mm}
{\rm GS}(J,p)
\subset \Bigg\{ p + \sum_{i \in \Lambda(J)}
                \alpha_i e_i \in M_a : \alpha_i \in \mathbb{R}
                ~\mbox{for all}~ i \in \Lambda(J)
        \Bigg\}.
\end{eqnarray*}

It remains to prove the reverse inclusion.
According to the structure of $J$ given in Definition
\ref{def:basic}, for each $i \in \Lambda(J)$, there exists a
real number $\beta_i \neq 0$ such that $p + \beta_i e_i \in J$.
In other words, for each $i \in \Lambda(J)$, there exists a
$u_i \in J$ such that $\beta_i e_i = u_i - p$.
Thus, if we assume that
\begin{eqnarray*}
p + \sum_{i \in \Lambda(J)} \alpha_i e_i \in M_a
\end{eqnarray*}
for some real numbers $\alpha_i$, then
\begin{eqnarray*}
p + \sum_{i \in \Lambda(J)} \alpha_i e_i
=   p + \sum_{i \in \Lambda(J)} \frac{\alpha_i}{\beta_i}
    (\beta_i e_i)
=   p + \sum_{i \in \Lambda(J)} \frac{\alpha_i}{\beta_i}
    (u_i - p)
\in {\rm GS}(J,p),
\end{eqnarray*}
since $u_i \in J$ for $i \in \Lambda(J)$, which implies that
\begin{eqnarray*}
\hspace*{3mm}
{\rm GS}(J,p)
\supset \Bigg\{ p + \sum_{i \in \Lambda(J)}
                \alpha_i e_i \in M_a : \alpha_i \in \mathbb{R}
                ~\mbox{for all}~ i \in \Lambda(J)
        \Bigg\}.
\end{eqnarray*}
We end the proof in this way.
\hfill$\Box$
\vspace*{5mm}

The following theorem may seem obvious, but it was derived in
an effort to achieve a similar result to the previous theorem
for more general sets.
We remark that a subset of $M_a$ with nonempty index set
contains at least two elements.


\begin{thm} \label{thm:3.7}
Assume that a bounded subset $E$ of $M_a$ contains at least
two elements and $p \in E$.
Then
\begin{eqnarray*}
{\rm GS}^2(E,p)
= \Bigg\{ p + \sum_{i \in \Lambda({\rm GS}^2(E,p))}
          \alpha_i e_i \in M_a :
          \alpha_i \in \mathbb{R} ~\mbox{for all}~
          i \in \Lambda \big( {\rm GS}^2(E,p) \big)
  \Bigg\}.
\end{eqnarray*}
\end{thm}

\noindent
\emph{Proof}.
According to Lemma \ref{lem:setahob}, if
$i \in \Lambda({\rm GS}^2(E,p))$, then
$\alpha_i e_i \in {\rm GS}^2(E,p) - p$ for all
$\alpha_i \in \mathbb{R}$.
Since ${\rm GS}^2(E,p) - p$ is a real vector space, if we set
$\Lambda_n = \{ i \in \Lambda({\rm GS}^2(E,p)) : i < n \}$,
then we have
\begin{eqnarray*}
\sum_{i \in \Lambda_n} \alpha_i e_i \in {\rm GS}^2(E,p) - p
\end{eqnarray*}
for all $n \in \mathbb{N}$ and all $\alpha_i \in \mathbb{R}$.
For now, with all $\alpha_i$ fixed, we define
$x_n = p + \sum\limits_{i \in \Lambda_n} \alpha_i e_i$ for any
$n \in \mathbb{N}$.
Then $\{ x_n \}$ is a sequence in ${\rm GS}^2(E,p)$.
When $\{ x_n \}$ converges in $M_a$, it holds that
\begin{eqnarray*}
p + \sum\limits_{i \in \Lambda({\rm GS}^2(E,p))} \alpha_i e_i
=   \lim_{n \to \infty} x_n
\in {\rm GS}^2(E,p),
\end{eqnarray*}
because ${\rm GS}^2(E,p)$ is closed by Lemma \ref{lem:3.8}
$(iii)$.
That is,
\begin{eqnarray*}
\Bigg\{ p + \sum_{i \in \Lambda({\rm GS}^2(E,p))}
        \alpha_i e_i \in M_a :
        \alpha_i \in \mathbb{R} ~\mbox{for all}~
        i \in \Lambda \big( {\rm GS}^2(E,p) \big)
\Bigg\} \subset {\rm GS}^2(E,p).
\end{eqnarray*}

It remains to prove the reverse inclusion.
Let $r$ be a positive real number with $E \subset B_r(p)$.
Assume that $x \in {\rm GS}^2(E,p)$ and note that
$\{ \frac{1}{a_i} e_i \}$ is a complete orthonormal sequence
in $M_a$ and
$\langle x-p, \frac{1}{a_i} e_i \rangle_a = a_i (x_i - p_i)$
for all $i$.
By Lemma \ref{lem:doolset} with ${\rm GS}^2(E,p) \cap B_r(p)$
in place of $E$ and considering Lemma \ref{lem:3.8} $(iv)$ and
Remark \ref{rem:new3.1}, we have
\begin{eqnarray} \label{eq:2021-4-20-1}
\begin{split}
x - p & = \sum_{i=1}^\infty
          \bigg\langle x-p, \frac{1}{a_i} e_i \bigg\rangle_a
          \frac{1}{a_i} e_i \\
      & = \sum_{i=1}^\infty a_i (x_i - p_i) \frac{1}{a_i} e_i \\
      & = \sum_{i \in \Lambda({\rm GS}^2(E,p))}
          a_i (x_i - p_i) \frac{1}{a_i} e_i \\
      & = \sum_{i \in \Lambda({\rm GS}^2(E,p))}
          \bigg\langle x-p, \frac{1}{a_i} e_i \bigg\rangle_a
          \frac{1}{a_i} e_i.
\end{split}
\end{eqnarray}
Since $x \in M_a$, it follows from Remark \ref{rem:new2.1} and
(\ref{eq:2021-4-20-1}) that
\begin{eqnarray} \label{eq:2021-4-20-2}
p + \sum_{i \in \Lambda({\rm GS}^2(E,p))}
\bigg\langle x-p, \frac{1}{a_i} e_i \bigg\rangle_a
\frac{1}{a_i} e_i = x \in M_a.
\end{eqnarray}

Finally, it follows from (\ref{eq:2021-4-20-1}) and
(\ref{eq:2021-4-20-2}) that
\begin{eqnarray*}
\begin{split}
x & =   p + \sum_{i \in \Lambda({\rm GS}^2(E,p))}
        \bigg\langle x-p, \frac{1}{a_i} e_i \bigg\rangle_a
        \frac{1}{a_i} e_i ~(\in M_a) \\
  & \in \Bigg\{ p + \sum_{i \in \Lambda({\rm GS}^2(E,p))}
                \alpha_i e_i \in M_a :
                \alpha_i \in \mathbb{R} ~\mbox{for all}~
                i \in \Lambda \big( {\rm GS}^2(E,p) \big)
        \Bigg\}
\end{split}
\end{eqnarray*}
for all $x \in {\rm GS}^2(E,p)$.
Thus, we can conclude that
\begin{eqnarray*}
\hspace*{3mm}
{\rm GS}^2(E,p)
\subset \Bigg\{ p + \sum_{i \in \Lambda({\rm GS}^2(E,p))}
                \alpha_i e_i \in M_a :
                \alpha_i \in \mathbb{R} ~\mbox{for all}~
                i \in \Lambda \big( {\rm GS}^2(E,p) \big)
        \Bigg\},
\end{eqnarray*}
which completes the proof.
\hfill$\Box$
\vspace*{5mm}

In the following theorem, we introduce an interesting inclusion
property of the second-order generalized linear span.


\begin{thm} \label{thm:setten}
Assume that a bounded subset $E$ of $M_a$ contains at least
two elements and $p \in E$.
Let $J$ be either a translation of a basic cylinder or a
translation of a basic interval that possesses the following
properties:
\begin{itemize}
\item[$(i)$]  There exists a $p' \in J \cap E$;
\item[$(ii)$] $\Lambda(J) \subset \Lambda({\rm GS}^2(E,p))$.
\end{itemize}
Then $J \subset {\rm GS}(J,p') \subset {\rm GS}^2(E,p)$.
\end{thm}

\noindent
\emph{Proof}.
It is obvious that $J \subset {\rm GS}(J,p')$.
Now we assume that $x \in {\rm GS}(J,p')$.
Then, according to Theorem \ref{thm:set9}, there exist some
real numbers $\alpha_i$ that satisfy
\begin{eqnarray} \label{eq:2021-10-3-1}
x = p' + \sum_{i \in \Lambda(J)} \alpha_i e_i \in M_a.
\end{eqnarray}
Since $p' \in E \subset {\rm GS}(E,p)$, it follows from Lemma
\ref{lem:setahob} and $(ii)$ that
$p' + \alpha_i e_i \in {\rm GS}^2(E,p)$ for all
$i \in \Lambda(J)$, \textit{i.e.},
\begin{eqnarray} \label{eq:2021-10-3-2}
\alpha_i e_i \in {\rm GS}^2(E,p) - p'
\end{eqnarray}
for all $i \in \Lambda(J)$.

Now we define
\begin{eqnarray*}
x_n := p' + \sum_{i \in \Lambda_n} \alpha_i e_i
\end{eqnarray*}
for any $n \in \mathbb{N}$, where
$\Lambda_n = \{ i \in \Lambda(J) : i < n \}$.
Since ${\rm GS}^2(E,p) - p'$ is a real vector space by Lemma
\ref{lem:3.8} $(i)$, it follows from (\ref{eq:2021-10-3-2})
that
\begin{eqnarray*}
\sum_{i \in \Lambda_n} \alpha_i e_i \in {\rm GS}^2(E,p) - p'
\end{eqnarray*}
for each $n \in \mathbb{N}$.
In other words,
\begin{eqnarray*}
x_n =   p' + \sum_{i \in \Lambda_n} \alpha_i e_i
    \in {\rm GS}^2(E,p)
\end{eqnarray*}
for each $n \in \mathbb{N}$.
Since ${\rm GS}^2(E,p)$ is closed by Lemma \ref{lem:3.8}
$(iii)$, it follows from (\ref{eq:2021-10-3-1}) that
\begin{eqnarray*}
x = \lim_{n \to \infty} x_n \in {\rm GS}^2(E,p),
\end{eqnarray*}
which implies that
$J \subset {\rm GS}(J,p') \subset {\rm GS}^2(E,p)$.
\hfill$\Box$
\vspace*{5mm}

Since in some ways index sets have some properties of dimensions
in vector space, the following theorem may seem obvious.


\begin{thm} \label{thm:setsib}
Assume that a bounded subset $E$ of $M_a$ contains at least
two elements and $p \in E$.
Then, $\Lambda({\rm GS}^2(E,p)) = \mathbb{N}$ if and only if\/
${\rm GS}^2(E,p) = M_a$.
\end{thm}

\noindent
\emph{Proof}.
Let $x$ be an arbitrary element of $M_a$.
There exist some real numbers $\alpha_i$ such that
\begin{eqnarray} \label{eq:2021-10-9}
x = \sum\limits_{i=1}^\infty \alpha_i e_i \in M_a.
\end{eqnarray}

If $\Lambda({\rm GS}^2(E,p)) = \mathbb{N}$, then it follows
from Lemma \ref{lem:setahob} that
\begin{eqnarray*}
p + \alpha_i e_i \in {\rm GS}^2(E,p)
\end{eqnarray*}
for all $i \in \mathbb{N}$.
In other words,
\begin{eqnarray*}
\alpha_i e_i \in {\rm GS}^2(E,p) - p
\end{eqnarray*}
for all $i \in \mathbb{N}$.

By Lemma \ref{lem:3.8} $(i)$, we get
\begin{eqnarray*}
x_n := \sum_{i=1}^n \alpha_i e_i \in {\rm GS}^2(E,p) - p
\end{eqnarray*}
for any $n \in \mathbb{N}$.
Due to Lemma \ref{lem:3.8} $(iii)$ and (\ref{eq:2021-10-9}),
we further obtain
\begin{eqnarray*}
x =   \sum_{i=1}^\infty \alpha_i e_i
  =   \lim_{n \to \infty} x_n
  \in {\rm GS}^2(E,p) - p,
\end{eqnarray*}
which implies that $M_a \subset {\rm GS}^2(E,p) - p$, or
equivalently, $M_a \subset {\rm GS}^2(E,p)$.

The reverse inclusion is trivial.
\hfill$\Box$


\section{Second-order extension of isometries}
\label{sec:7}

It was proved in Theorem \ref{thm:new3} that the domain of
$d_a$-isometry $f : E_1 \to E_2$ can be extended to the
first-order generalized linear span ${\rm GS}(E_1,p)$ whenever
$E_1$ is a nonempty bounded subset of $M_a$, whether degenerate
or non-degenerate.

Now we generalize Theorem \ref{thm:new3} in the following
theorem.
More precisely, we prove that the domain of $f$ can be extended
to its second-order generalized linear span ${\rm GS}^2(E_1,p)$.
It follows from Lemma \ref{lem:3.8} $(iii)$ that
${\rm GS}^2(E_1,p) = \overline{{\rm GS}(E_1,p)}$.
Therefore, the following theorem is a further generalization
of \cite[Theorem 2.2]{jung1}.
In the proof, we use the fact that ${\rm GS}^n(E_1,p) - p$ is
a real vector space.


\begin{thm} \label{thm:setahob}
Let $E_1$ be a bounded subset of $M_a$ that is $d_a$-isometric to
a subset $E_2$ of $M_a$ via a surjective $d_a$-isometry
$f : E_1 \to E_2$.
Assume that $p$ and $q$ are elements of $E_1$ and $E_2$, which
satisfy $q = f(p)$.
The function $F_2 : {\rm GS}^2(E_1,p) \to M_a$ is a $d_a$-isometry and the function
$T_{-q} \circ F_2 \circ T_p : {\rm GS}^2(E_1,p) - p \to M_a$
is linear.
In particular, $F_2$ is an extension of $F$.
\end{thm}

\noindent
\emph{Proof}.
$(a)$ Suppose $r$ is a positive real number satisfying
$E_1 \subset B_r(p)$.
Referring to the changes presented in the table below
and following the first part of proof of Theorem \ref{thm:new3},
we can easily prove that $F_2$ is a $d_a$-isometry.

\begin{center}
\begin{tabular}{ccccccc}
\hline\hline
Theorem \ref{thm:new3}: & $E_1$                         &
${\rm GS}(E_1,p)$   & $f$ & $F$   &
Definition \ref{def:ext}        & Lemma \ref{lem:new2} \\
\hline
Here:                   & ${\rm GS}(E_1,p) \cap B_r(p)$ &
${\rm GS}^2(E_1,p)$ & $F$ & $F_2$ &
Definition \ref{def:6.1} & Lemma \ref{lem:6.1}  \\
\hline\hline
\end{tabular}
\end{center}

$(b)$ We prove the linearity of
$T_{-q} \circ F_n \circ T_p : {\rm GS}^n(E_1,p) - p \to M_a$
in the more general setting for $n \geq 2$.
Referring to the changes presented in the table below
and following $(d)$ of the proof of Theorem \ref{thm:new3},
we can easily prove that $T_{-q} \circ F_n \circ T_p$ is linear.

\begin{center}
\begin{tabular}{cccc}
\hline\hline
Theorem \ref{thm:new3}: & ${\rm GS}(E_1,p)$   & $F$
& (\ref{eq:2018-9-28}) \\
\hline
Here:                   & ${\rm GS}^n(E_1,p)$ & $F_n$
& Lemma \ref{lem:6.1}  \\
\hline\hline
\end{tabular}
\end{center}
\vspace*{2mm}

$(c)$ According to Definition \ref{def:6.1} $(i)$, for any
$m \in \mathbb{N}$, $x_{ij} \in {\rm GS}(E_1,p) \cap B_r(p)$,
and any $\alpha_{ij} \in \mathbb{R}$ with
$\sum\limits_{i=1}^m \sum\limits_{j=1}^\infty \alpha_{ij}
 (x_{ij} - p) \in M_a$, there exists a
$u \in {\rm GS}^2(E_1,p)$ satisfying
\begin{eqnarray} \label{eq:3.6a}
u - p = \sum_{i=1}^m \sum_{j=1}^\infty
        \alpha_{ij} (x_{ij} - p) \in M_a.
\end{eqnarray}
Due to Definition \ref{def:6.1} $(ii)$, we further have
\begin{eqnarray} \label{eq:3.7a}
(T_{-q} \circ F_2 \circ T_p)(u - p)
= \sum_{i=1}^m \sum_{j=1}^\infty
  \alpha_{ij} (T_{-q} \circ F \circ T_p)(x_{ij} - p).
\end{eqnarray}

If we set $\alpha_{11} = 1$, $\alpha_{ij} = 0$ for each
$(i,j) \neq (1,1)$, and $x_{11} = x$ in (\ref{eq:3.6a}) and
(\ref{eq:3.7a}) to see
\begin{eqnarray} \label{eq:2021-9-26}
(T_{-q} \circ F_2 \circ T_p)(x-p)
= (T_{-q} \circ F \circ T_p)(x-p)
\end{eqnarray}
for all $x \in {\rm GS}(E_1,p) \cap B_r(p)$.

Let $w$ be an arbitrary element of ${\rm GS}(E_1,p)$.
Then, $w - p \in {\rm GS}(E_1,p) - p$.
Since ${\rm GS}(E_1,p) - p$ is a real vector space and
$B_r(p) - p = B_r(0)$, there exists a (sufficiently small)
real number $\mu \neq 0$ such that
\begin{eqnarray*}
\mu (w - p) \in ({\rm GS}(E_1,p) - p) \cap (B_r(p) - p).
\end{eqnarray*}
Hence, by (\ref{eq:2021-11-2}), we can choose a
$v \in {\rm GS}(E_1,p) \cap B_r(p)$ such that
$\mu (w - p) = v - p$.
Since both $T_{-q} \circ F_2 \circ T_p$ and
$T_{-q} \circ F \circ T_p$ are linear and
${\rm GS}(E_1,p) \subset {\rm GS}^2(E_1,p)$, it follows from
(\ref{eq:2021-9-26}) that
\begin{eqnarray*}
\begin{split}
\mu (T_{-q} \circ F_2 \circ T_p)(w - p)
& = (T_{-q} \circ F_2 \circ T_p)(\mu (w - p)) \\
& = (T_{-q} \circ F_2 \circ T_p)(v - p) \\
& = (T_{-q} \circ F \circ T_p)(v - p) \\
& = (T_{-q} \circ F \circ T_p)(\mu (w - p)) \\
& = \mu (T_{-q} \circ F \circ T_p)(w - p).
\end{split}
\end{eqnarray*}
Therefore, it follows that
$(T_{-q} \circ F_2 \circ T_p)(w-p)
 = (T_{-q} \circ F \circ T_p)(w-p)$ for all
$w \in {\rm GS}(E_1,p)$, \textit{i.e.}, $F_2(w) = F(w)$ for
all $w \in {\rm GS}(E_1,p)$.
In other words, $F_2$ is an extension of $F$.
Also, because of Theorem \ref{thm:new3}, we see that $F_2$ is
obviously an extension of $f$.
\hfill$\Box$
\vspace*{5mm}

On account of Remark \ref{rem:new3.1}, it holds that
\begin{eqnarray*}
{\rm GS}^2(E_1,p) = \cdots = {\rm GS}^{n-1}(E_1,p)
= {\rm GS}^n(E_1,p)
\end{eqnarray*}
for every integer $n \geq 3$.
According to this formula, the assertion of the following
theorem seems obvious, but since the proof is not long, we
introduce the proof here.


\begin{thm} \label{thm:setyol}
Let $E_1$ be a bounded subset of $M_a$ that is $d_a$-isometric to
a subset $E_2$ of $M_a$ via a surjective $d_a$-isometry
$f : E_1 \to E_2$.
Assume that $p$ and $q$ are elements of $E_1$ and $E_2$, which
satisfy $q = f(p)$.
Then $F_n$ is identically the same as $F_2$ for any integer
$n \geq 3$, where $F_2$ and $F_n$ are defined in Definition
\ref{def:6.1}.
\end{thm}

\noindent
\emph{Proof}.
Let $r$ be a fixed positive real number satisfying
$E_1 \subset B_r(p)$.
We assume that $F_2 \equiv F_3 \equiv \cdots \equiv F_{n-1}$
on ${\rm GS}^2(E_1,p)$.
Let $x$ be an arbitrary element of ${\rm GS}^n(E_1,p)$.
Then, in view of (\ref{eq:2021-11-2}), there exist a real
number $\mu \neq 0$ and an element $u$ of
${\rm GS}^n(E_1,p) \cap B_r(p)$ such that
\begin{eqnarray*}
u - p = \mu (x - p)
\in ({\rm GS}^n(E_1,p) - p) \cap (B_r(p) - p).
\end{eqnarray*}
If we put $\alpha_{11} = 1$, $\alpha_{ij} = 0$ for all
$(i,j) \neq (1,1)$, and $x_{11} = v$ in Definition \ref{def:6.1}
$(ii)$, then we get
\begin{eqnarray} \label{eq:2021-10-29}
(T_{-q} \circ F_n \circ T_p)(v - p)
= (T_{-q} \circ F_{n-1} \circ T_p)(v - p)
\end{eqnarray}
for all
$v \in {\rm GS}^{n-1}(E_1,p) \cap B_r(p)
   =   {\rm GS}^n(E_1,p) \cap B_r(p)$ by Remark
\ref{rem:new3.1}.

Since $T_{-q} \circ F_n \circ T_p$ is linear by $(b)$ in the
proof of Theorem \ref{thm:setahob}, it follows from
(\ref{eq:2021-10-29}) and our assumption that
\begin{eqnarray*}
\begin{split}
\mu (T_{-q} \circ F_n \circ T_p)(x - p)
& = (T_{-q} \circ F_n \circ T_p)(u - p) \\
& = (T_{-q} \circ F_{n-1} \circ T_p)(u - p) \\
& = (T_{-q} \circ F_2 \circ T_p)(u - p) \\
& = \mu (T_{-q} \circ F_2 \circ T_p)(x - p),
\end{split}
\end{eqnarray*}
\textit{i.e.}, $F_n(x) = F_2(x)$ for every
$x \in {\rm GS}^n(E_1,p) = {\rm GS}^2(E_1,p)$.
By mathematical induction, we conclude that $F_n$ is
identically the same as $F_2$ for every integer $n \geq 3$.
\hfill$\Box$
\vspace*{5mm}

We note that $\{ \frac{1}{a_i} e_i \}$ is a complete
orthonormal sequence in $M_a$.
Assume that $J$ is either a translation of a basic cylinder or
a translation of a basic interval, and $p$ is an element of
$J$.
On account of Theorem \ref{thm:set9}, we notice that
$\Lambda(J) = \Lambda({\rm GS}(J,p))$.
On the other hand, if we apply Lemma \ref{lem:doolset} with
${\rm GS}(J,p) \cap B_r(p)$ in place of $E$, where $r$ is a
positive real number with $J \subset B_r(p)$, then
$i \not\in \Lambda({\rm GS}(J,p) \cap B_r(p))$ if and only if
$u_i = p_i$ for all
$u \in {\rm GS}({\rm GS}(J,p) \cap B_r(p), p)$.
Equivalently, it holds that
$i \not\in \Lambda({\rm GS}(J,p))$ if and only if $u_i = p_i$
for all $u \in {\rm GS}^2(J,p)$.
Hence, we have
\begin{eqnarray} \label{eq:2021-10-21-2}
\begin{split}
u - p
& = \sum_{i=1}^\infty
    \left\langle u - p,\, \frac{1}{a_i} e_i \right\rangle_a
    \frac{1}{a_i} e_i \\
& = \sum_{i=1}^\infty a_i (u_i - p_i) \frac{1}{a_i} e_i \\
& = \sum_{i \in \Lambda({\rm GS}(J,p))}
    a_i (u_i - p_i) \frac{1}{a_i} e_i \\
& = \sum_{i \in \Lambda({\rm GS}(J,p))}
    \left\langle u - p,\, \frac{1}{a_i} e_i \right\rangle_a
    \frac{1}{a_i} e_i \\
& = \sum_{i \in \Lambda(J)}
    \left\langle u - p,\, \frac{1}{a_i} e_i \right\rangle_a
    \frac{1}{a_i} e_i
\end{split}
\end{eqnarray}
for all $u \in {\rm GS}^2(J,p) = {\rm GS}^n(J,p)$, where
$n \geq 2$.

Using a similar approach to the proof of
\cite[Theorem 2.4]{jungkim}, we can apply Lemma \ref{lem:6.1}
to prove the following theorem.


\begin{thm} \label{thm:new4}
Assume that $J$ is either a translation of a basic cylinder or
a translation of a basic interval, $K$ is a subset of $M_a$,
and that there exists a surjective $d_a$-isometry $f : J \to K$.
Suppose $p$ is an element of $J$ and $q$ is an element of $K$
with $q = f(p)$.
For any $n \in \mathbb{N}$, the $d_a$-isometry
$F_n : {\rm GS}^n(J,p) \to M_a$ given in Definition
\ref{def:6.1} satisfies
\begin{eqnarray*}
(T_{-q} \circ F_n \circ T_p)(u-p)
= \sum_{i \in \Lambda(J)}
  \bigg\langle u-p, \frac{1}{a_i} e_i \bigg\rangle_a
  \frac{1}{a_i} (T_{-q} \circ F_n \circ T_p)(e_i)
\end{eqnarray*}
for all $u \in {\rm GS}^n(J,p)$.
\end{thm}

\noindent
\emph{Proof}.
Since $p + e_i \in {\rm GS}^n(J,p)$ for each $i \in \Lambda(J)$,
it follows from Lemma \ref{lem:6.1} that
\begin{eqnarray}
\begin{split}
& \Bigg\langle (T_{-q} \circ F_n \circ T_p)(u-p) -
               \sum_{i \in \Lambda(J)}
               \left\langle u-p, \frac{1}{a_i} e_i \right\rangle_a
               \frac{1}{a_i}
               (T_{-q} \circ F_n \circ T_p)(e_i), \\
& \hspace*{3mm}
    \left.     (T_{-q} \circ F_n \circ T_p)(u-p) -
               \sum_{j \in \Lambda(J)}
               \left\langle u-p, \frac{1}{a_j} e_j \right\rangle_a
               \frac{1}{a_j}
               (T_{-q} \circ F_n \circ T_p)(e_j)
  \right\rangle_a \\
& \hspace*{3mm}
  = \big\langle (T_{-q} \circ F_n \circ T_p)(u-p),\,
                (T_{-q} \circ F_n \circ T_p)(u-p)
    \big\rangle_a \\
& \hspace*{8mm}
    -\,\sum_{j \in \Lambda(J)}
    \left\langle u-p, \frac{1}{a_j} e_j \right\rangle_a
    \frac{1}{a_j}
    \big\langle (T_{-q} \circ F_n \circ T_p)(u-p),\,
                (T_{-q} \circ F_n \circ T_p)(e_j)
    \big\rangle_a \\
& \hspace*{8mm}
    -\,\sum_{i \in \Lambda(J)}
    \left\langle u-p, \frac{1}{a_i} e_i \right\rangle_a
    \frac{1}{a_i}
    \big\langle (T_{-q} \circ F_n \circ T_p)(e_i),\,
                (T_{-q} \circ F_n \circ T_p)(u-p)
    \big\rangle_a \\
& \hspace*{8mm}
    +\,\sum_{i \in \Lambda(J)} \sum_{j \in \Lambda(J)}
    \left\langle u-p, \frac{1}{a_i} e_i \right\rangle_a
    \left\langle u-p, \frac{1}{a_j} e_j \right\rangle_a
    \times
    \label{eq:3.3} \\
& \hspace*{34mm} \times
    \frac{1}{a_i a_j}
    \big\langle (T_{-q} \circ F_n \circ T_p)(e_i),\,
                (T_{-q} \circ F_n \circ T_p)(e_j)
    \big\rangle_a \\
& \hspace*{3mm}
  = \big\langle u-p, u-p \big\rangle_a -
    \sum_{j \in \Lambda(J)}
    \left\langle u-p, \frac{1}{a_j} e_j \right\rangle_a
    \left\langle u-p, \frac{1}{a_j} e_j \right\rangle_a \\
& \hspace*{8mm}
    -\,\sum_{i \in \Lambda(J)}
    \left\langle u-p, \frac{1}{a_i} e_i \right\rangle_a
    \left\langle \frac{1}{a_i} e_i, u-p \right\rangle_a \\
& \hspace*{8mm}
    +\,\sum_{i \in \Lambda(J)} \sum_{j \in \Lambda(J)}
    \left\langle u-p, \frac{1}{a_i} e_i \right\rangle_a
    \left\langle u-p, \frac{1}{a_j} e_j \right\rangle_a
    \left\langle \frac{1}{a_i} e_i,\, \frac{1}{a_j} e_j
    \right\rangle_a \\
& \hspace*{3mm}
  = \big\langle u-p, u-p \big\rangle_a -
    \sum_{j \in \Lambda(J)}
    \left\langle u-p, \frac{1}{a_j} e_j \right\rangle_a
    \left\langle u-p, \frac{1}{a_j} e_j \right\rangle_a
\end{split}
\end{eqnarray}
for all $u \in {\rm GS}^n(J,p)$, since
$\big\{ \frac{1}{a_i} e_i \big\}$ is an orthonormal sequence
in $M_a$.

Furthermore, we note that each $u \in {\rm GS}^n(J,p)$ has the
expression given in (\ref{eq:2021-10-21-2}).
Hence, if we replace $u-p$ in (\ref{eq:3.3}) with the expression
(\ref{eq:2021-10-21-2}), then we have
\begin{eqnarray*}
\Bigg\| (T_{-q} \circ F_n \circ T_p)(u-p) -
        \sum_{i \in \Lambda(J)}
        \left\langle u-p, \frac{1}{a_i} e_i \right\rangle_a
        \frac{1}{a_i} (T_{-q} \circ F_n \circ T_p)(e_i)
\Bigg\|_a^2 = 0
\end{eqnarray*}
for all $u \in {\rm GS}^n(J,p)$, which implies the validity of our
assertion.
\hfill$\Box$
\vspace*{5mm}

According to the following theorem, the image of the
first-order generalized linear span of $E_1$ with respect to
$p$ under the $d_a$-isometry $F$ is just the first-order generalized
linear span of $F(E_1)$ with respect to $F(p)$.
This assertion holds also for the second-order generalized
linear span and $F_2$.
According to Remark \ref{rem:new3.1} and Theorem
\ref{thm:setyol}, the argument of the following theorem only
makes sense when $n = 1$ or $2$.


\begin{thm} \label{thm:2018-11-7}
Assume that $E_1$ and $E_2$ are bounded subsets of $M_a$ that
are $d_a$-isometric to each other via a surjective
$d_a$-isometry $f : E_1 \to E_2$.
Suppose $p$ is an element of $E_1$ and $q$ is an element of
$E_2$ with $q = f(p)$.
If $F_n : {\rm GS}^n(E_1,p) \to M_a$ is the extension of $f$
defined in Definition \ref{def:6.1}, then
${\rm GS}^n(E_2,q) = F_n({\rm GS}^n(E_1,p))$ for every
$n \in \mathbb{N}$.
\end{thm}

\noindent
\emph{Proof}.
$(a)$
First, we prove that our assertion is true for $n = 1$,
\textit{i.e.}, we prove that
${\rm GS}(E_2,q) = F({\rm GS}(E_1,p))$.
Let $r$ be a fixed positive real number satisfying
$E_1 \subset B_r(p)$.

$(b)$
Due to Definition \ref{def:gls}, for any
$y \in F({\rm GS}(E_1,p))$, there exists an element
$x \in {\rm GS}(E_1,p)$ with
\begin{eqnarray*}
y = F(x)
  = F \bigg( p + \sum_{i=1}^m \sum_{j=1}^\infty
             \alpha_{ij} \big( u_{ij} - p \big) \bigg)
\end{eqnarray*}
for some $m \in \mathbb{N}$, $u_{ij} \in E_1 \cap B_r(p)$, and
some $\alpha_{ij} \in \mathbb{R}$ with
$x = p + \sum\limits_{i=1}^m \sum\limits_{j=1}^\infty
 \alpha_{ij} (u_{ij} - p) \in M_a$.

On the other hand, by Definition \ref{def:ext}, we have
\begin{eqnarray*}
(T_{-q} \circ F \circ T_p)
\bigg( \sum_{i=1}^m \sum_{j=1}^\infty
       \alpha_{ij} \big( u_{ij} - p \big)
\bigg)
= \sum_{i=1}^m \sum_{j=1}^\infty
  \alpha_{ij} (T_{-q} \circ f \circ T_p) (u_{ij} - p)
\end{eqnarray*}
which is equivalent to
\begin{eqnarray*}
F(x) - q
= F \bigg( p + \sum_{i=1}^m \sum_{j=1}^\infty
           \alpha_{ij} \big( u_{ij} - p \big)
    \bigg) - q
= \sum_{i=1}^m \sum_{j=1}^\infty
  \alpha_{ij} \big( f(u_{ij}) - q \big).
\end{eqnarray*}
Since $u_{ij} \in E_1 = E_1 \cap B_r(p)$ for all $i$ and $j$,
it holds that $f(u_{ij}) \in f(E_1) = E_2$ for each $i$ and
$j$.
Moreover, since $u_{ij} \in E_1 \cap B_r(p)$ for all $i$ and
$j$, it follows from Lemma \ref{lem:new2} that
\begin{eqnarray*}
\begin{split}
\| f(u_{ij}) - q \|_a^2
& = \| (T_{-q} \circ f \circ T_p)(u_{ij} - p) \|_a^2 \\
& = \big\langle (T_{-q} \circ f \circ T_p)(u_{ij} - p),\,
                (T_{-q} \circ f \circ T_p)(u_{ij} - p)
    \big\rangle_a \\
& = \langle u_{ij} - p, u_{ij} - p \rangle_a \\
& = \| u_{ij} - p \|_a^2 \\
& < r^2
\end{split}
\end{eqnarray*}
for all $i$ and $j$.
Hence, $f(u_{ij}) \in E_2 \cap B_r(q)$ for all $i$ and $j$.

Furthermore, it follows from Lemma \ref{lem:new2} that
\begin{eqnarray*}
\begin{split}
& \left\| \sum_{i=1}^m \sum_{j=1}^\infty
          \alpha_{ij} \big( f(u_{ij}) - q \big)
  \right\|_a^2 \\
& \hspace*{3mm}
  = \left\| \sum_{i=1}^m \sum_{j=1}^\infty
            \alpha_{ij} (T_{-q} \circ f \circ T_p)(u_{ij} - p)
    \right\|_a^2 \\
& \hspace*{3mm}
  = \left\langle
    \sum_{i=1}^m \sum_{j=1}^\infty
    \alpha_{ij} (T_{-q} \circ f \circ T_p)(u_{ij} - p),\,
    \sum_{k=1}^m \sum_{\ell=1}^\infty
    \alpha_{k\ell} (T_{-q} \circ f \circ T_p)(u_{k\ell} - p)
    \right\rangle_a \\
& \hspace*{3mm}
  = \sum_{i=1}^m \sum_{k=1}^m
    \sum_{j=1}^\infty \alpha_{ij}
    \sum_{\ell=1}^\infty \alpha_{k\ell}
    \big\langle (T_{-q} \circ f \circ T_p)(u_{ij} - p),\,
                (T_{-q} \circ f \circ T_p)(u_{k\ell} - p)
    \big\rangle_a \\
& \hspace*{3mm}
  = \sum_{i=1}^m \sum_{k=1}^m
    \sum_{j=1}^\infty \alpha_{ij}
    \sum_{\ell=1}^\infty \alpha_{k\ell}
    \langle u_{ij} - p, u_{k\ell} - p \rangle_a \\
& \hspace*{3mm}
  = \left\langle
    \sum_{i=1}^m \sum_{j=1}^\infty \alpha_{ij} (u_{ij} - p),\,
    \sum_{k=1}^m \sum_{\ell=1}^\infty
    \alpha_{k\ell} (u_{k\ell} - p)
    \right\rangle_a \\
& \hspace*{3mm}
  = \left\|
    \sum_{i=1}^m \sum_{j=1}^\infty \alpha_{ij} (u_{ij} - p)
    \right\|_a^2 \\
& \hspace*{3mm}
  < \infty,
\end{split}
\end{eqnarray*}
since
$\sum\limits_{i=1}^m \sum\limits_{j=1}^\infty
 \alpha_{ij} (u_{ij} - p) = x - p \in M_a$.
Thus, on account of Remark \ref{rem:new2.1}, we see that
$\sum\limits_{i=1}^m \sum\limits_{j=1}^\infty
 \alpha_{ij} (f(u_{ij}) - q) \in M_a$.
Therefore, in view of Definition \ref{def:gls}, we get
\begin{eqnarray*}
y =   F(x)
  =   q + \sum_{i=1}^m \sum_{j=1}^\infty
      \alpha_{ij} \big( f(u_{ij}) - q \big)
  \in {\rm GS}(E_2,q)
\end{eqnarray*}
and we conclude that
$F({\rm GS}(E_1,p)) \subset {\rm GS}(E_2,q)$.

$(c)$ Now we assume that $y \in {\rm GS}(E_2,q)$.
By Definition \ref{def:gls}, there exist some
$m \in \mathbb{N}$, $v_{ij} \in E_2 \cap B_r(q)$, and some
$\alpha_{ij} \in \mathbb{R}$ such that
$y - q = \sum\limits_{i=1}^m \sum\limits_{j=1}^\infty
 \alpha_{ij} (v_{ij} - q) \in M_a$.
Since $f : E_1 \to E_2$ is surjective, there exists a
$u_{ij} \in E_1$ satisfying $v_{ij} = f(u_{ij})$ for any $i$
and $j$.
Moreover, by Lemma \ref{lem:new2}, we have
\begin{eqnarray*}
\begin{split}
\| u_{ij} - p \|_a^2
& = \langle u_{ij} - p, u_{ij} - p \rangle_a \\
& = \big\langle (T_{-q} \circ f \circ T_p)(u_{ij} - p),\,
                (T_{-q} \circ f \circ T_p)(u_{ij} - p)
    \big\rangle_a \\
& = \big\langle f(u_{ij}) - q, f(u_{ij}) - q \big\rangle_a \\
& = \langle v_{ij} - q, v_{ij} - q \rangle_a \\
& = \| v_{ij} - q \|_a^2 \\
& < r^2
\end{split}
\end{eqnarray*}
for any $i$ and $j$.
So we conclude that $u_{ij} \in E_1 \cap B_r(p)$ and
$v_{ij} = f(u_{ij})$ for all $i$ and $j$.

On the other hand, using Lemma \ref{lem:new2}, we have
\begin{eqnarray*}
\begin{split}
& \bigg\|
  \sum_{i=1}^m \sum_{j=1}^\infty \alpha_{ij} (u_{ij} - p)
  \bigg\|_a^2 \\
& \hspace*{3mm}
  = \left\langle
    \sum_{i=1}^m \sum_{j=1}^\infty \alpha_{ij} (u_{ij} - p),\,
    \sum_{k=1}^m \sum_{\ell=1}^\infty
    \alpha_{k\ell} (u_{k\ell} - p)
    \right\rangle_a \\
& \hspace*{3mm}
  = \sum_{i=1}^m \sum_{k=1}^m
    \sum_{j=1}^\infty \alpha_{ij}
    \sum_{\ell=1}^\infty \alpha_{k\ell}
    \big\langle u_{ij} - p,\, u_{k\ell} - p \big\rangle_a \\
& \hspace*{3mm}
  = \sum_{i=1}^m \sum_{k=1}^m
    \sum_{j=1}^\infty \alpha_{ij}
    \sum_{\ell=1}^\infty \alpha_{k\ell}
    \big\langle (T_{-q} \circ f \circ T_p)(u_{ij} - p),\,
                (T_{-q} \circ f \circ T_p)(u_{k\ell} - p)
    \big\rangle_a \\
& \hspace*{3mm}
  = \left\langle
    \sum_{i=1}^m \sum_{j=1}^\infty \alpha_{ij}
    (T_{-q} \circ f \circ T_p)(u_{ij} - p),\,
    \sum_{k=1}^m \sum_{\ell=1}^\infty \alpha_{k\ell}
    (T_{-q} \circ f \circ T_p)(u_{k\ell} - p)
    \right\rangle_a \\
& \hspace*{3mm}
  = \bigg\|
    \sum_{i=1}^m \sum_{j=1}^\infty
    \alpha_{ij} (T_{-q} \circ f \circ T_p)(u_{ij} - p)
    \bigg\|_a^2 \\
& \hspace*{3mm}
  = \bigg\|
    \sum_{i=1}^m \sum_{j=1}^\infty
    \alpha_{ij} \big( f(u_{ij}) - q \big)
    \bigg\|_a^2 \\
& \hspace*{3mm}
  = \bigg\|
    \sum_{i=1}^m \sum_{j=1}^\infty \alpha_{ij} (v_{ij} - q)
    \bigg\|_a^2 \\
& \hspace*{3mm}
  < \infty,
\end{split}
\end{eqnarray*}
since
$\sum\limits_{i=1}^m \sum\limits_{j=1}^\infty
 \alpha_{ij} (v_{ij} - q) = y - q \in M_a$.
Thus, Remark \ref{rem:new2.1} implies that
$\sum\limits_{i=1}^m \sum\limits_{j=1}^\infty
 \alpha_{ij} (u_{ij} - p) \in M_a$.

Hence, it follows from Definition \ref{def:ext} that
\begin{eqnarray*}
\begin{split}
y & =   q + \sum_{i=1}^m \sum_{j=1}^\infty
        \alpha_{ij} \big( f(u_{ij}) - q \big) \\
  & =   q + \sum_{i=1}^m \sum_{j=1}^\infty
        \alpha_{ij} (T_{-q} \circ f \circ T_p)(u_{ij} - p) \\
  & =   q + (T_{-q} \circ F \circ T_p)
        \bigg(
        \sum_{i=1}^m \sum_{j=1}^\infty \alpha_{ij} (u_{ij} - p)
        \bigg) \\
  & =   F \bigg( p + \sum_{i=1}^m \sum_{j=1}^\infty
                 \alpha_{ij} (u_{ij} - p)
          \bigg) \\
  & \in F({\rm GS}(E_1,p)).
\end{split}
\end{eqnarray*}
Thus, we conclude that
${\rm GS}(E_2,q) \subset F({\rm GS}(E_1,p))$.

$(d)$
Similarly, referring to the changes presented in the tables
below and following the previous parts $(b)$ and $(c)$ in this
proof, we can prove that
${\rm GS}^2(E_2,q) = F_2({\rm GS}^2(E_1,p))$.

\begin{center}
\begin{tabular}{ccccccc}
\hline\hline
The case $n = 1$:   & $E_1$               & $E_2$             &
${\rm GS}(E_1,p)$   & ${\rm GS}(E_2,q)$   & $f$ & $F$   \\
\hline
The case $n = 2$:   & ${\rm GS}(E_1,p)$   & ${\rm GS}(E_2,q)$ &
${\rm GS}^2(E_1,p)$ & ${\rm GS}^2(E_2,q)$ & $F$ & $F_2$ \\
\hline\hline
\end{tabular}
\end{center}

\begin{center}
\begin{tabular}{cccc}
\hline\hline
The case $n = 1$:   & Definition \ref{def:gls}
& Definition \ref{def:ext}        & Lemma \ref{lem:new2}   \\
\hline
The case $n = 2$:   & Definition \ref{def:6.1} $(i)$
& Definition \ref{def:6.1} $(ii)$ & $(\ref{eq:2018-9-28})$ \\
\hline\hline
\end{tabular}
\end{center}

$(e)$
Finally, according to Remark \ref{rem:new3.1}, Theorem
\ref{thm:setyol}, and $(d)$, we further have
\begin{eqnarray*}
{\rm GS}^n(E_2,q) = {\rm GS}^2(E_2,q) = F_2({\rm GS}^2(E_1,p))
= F_n({\rm GS}^n(E_1,p))
\end{eqnarray*}
for any integer $n \geq 3$.
\hfill$\Box$


\section{Extension of isometries to entire space}
\label{sec:8}

Let $I^\omega = \prod\limits_{i=1}^\infty I$ be the
\emph{Hilbert cube}, where $I = [0,1]$ is the unit closed
interval.
From now on, we assume that $E_1$ and $E_2$ are nonempty
subsets of $I^\omega$.
They are bounded, of course.

In our main theorem (Theorem \ref{thm:3.17}), we will prove
that the domain of local $d_a$-isometry $f : E_1 \to E_2$ can
be extended to any real Hilbert space including the domain
$E_1$.


\begin{defn} \label{def:3.5}
Let $E_1$ be a nonempty subset of $I^\omega$ that is
$d_a$-isometric to a subset $E_2$ of $I^\omega$ via a
surjective $d_a$-isometry $f : E_1 \to E_2$.
Let $p$ be an element of $E_1$ and $q$ an element of $E_2$
with $q = f(p)$.
Assume that
$\mathbb{N} \!\setminus\! \Lambda({\rm GS}^2(E_1,p))$ is a
nonempty set and $\Lambda$ is a set with
$\Lambda({\rm GS}^2(E_1,p)) \subset \Lambda \subset \mathbb{N}$.
Suppose $\{ \frac{1}{a_i} e'_i \}_{i \in \Lambda}$ is an
orthonormal sequence and
$e'_i = (T_{-q} \circ F_2 \circ T_p)(e_i)$ for each
$i \in \Lambda({\rm GS}^2(E_1,p))$, where
$F_2 : {\rm GS}^2(E_1,p) \to M_a$ is defined in Definition
\ref{def:6.1}.
Let $p_i$ be the $i$th component of $p$.
We define a basic cylinder or a basic interval $\tilde{J}$ by
\begin{eqnarray*}
\tilde{J} = \prod_{i=1}^\infty\, \tilde{J}_i,
~~~\mbox{where}~~~
\tilde{J}_i = \left\{\begin{array}{ll}
                        [0, 1]
                        & (\mbox{for}~i \in \Lambda),
                        \vspace*{1mm}\\
                        \{ p_i \}
                        & (\mbox{for}~i \not\in \Lambda).
                     \end{array}
              \right.
\end{eqnarray*}
Moreover, we define the function
$G_2 : {\rm GS}^2(\tilde{J},p) \to M_a$ by
\begin{eqnarray} \label{eq:2021-10-21-1}
(T_{-q} \circ G_2 \circ T_p)(u - p)
= \sum_{i \in \Lambda(\tilde{J})}
  \left\langle u - p,\, \frac{1}{a_i} e_i \right\rangle_a
  \frac{1}{a_i} e'_i
\end{eqnarray}
for all $u \in {\rm GS}^2(\tilde{J},p)$.
\end{defn}

Since $\{ \frac{1}{a_i} e_i \}_{i \in \Lambda({\rm GS}^2(E_1,p))}$
is a complete orthonormal sequence in the Hilbert space
${\rm GS}^2(E_1,p) - p$ and
${\rm GS}^2(E_2,q) = F_2({\rm GS}^2(E_1,p))$, the following
lemma seems obvious, but its proof is not very simple, so we
introduce it here.


\begin{lem} \label{lem:8.1}
Assume that $E_1$ is a bounded subset of $I^\omega$ that
contains at least two elements and is $d_a$-isometric to a
subset $E_2$ of $I^\omega$ via a surjective $d_a$-isometry
$f : E_1 \to E_2$.
Moreover, assume that $p$ and $q$ are elements of $E_1$ and
$E_2$ that satisfy $q = f(p)$.
Then
$\{ \frac{1}{a_i} (T_{-q} \circ F_2 \circ T_p)(e_i) \}
 _{i \in \Lambda({\rm GS}^2(E_1,p))}$ is a complete orthonormal
sequence in ${\rm GS}^2(E_2,q) - q$, where $F_2$ is the
extension of $f$ defined in Definition \ref{def:6.1}.
\end{lem}

\noindent
\emph{Proof}.
It follows from Lemma \ref{lem:setahob} that
$p + e_i \in {\rm GS}^2(E_1,p)$ or
$e_i \in {\rm GS}^2(E_1,p) - p$ for any
$i \in \Lambda({\rm GS}^2(E_1,p))$.
Hence, for any $i \in \Lambda({\rm GS}^2(E_1,p))$, there exists
a $u_i \in {\rm GS}^2(E_1,p)$ such that $e_i = u_i - p$.
Moreover, we have
\begin{eqnarray*}
\| u_i - p \|_a^2
= \| e_i \|_a^2
= a_i^2 \left\langle \frac{1}{a_i} e_i, \frac{1}{a_i} e_i
        \right\rangle_a
= a_i^2,
\end{eqnarray*}
which implies that $u_i \in B_r(p)$, where we set
$r = \sup\limits_{i \in \Lambda({\rm GS}^2(E_1,p))} a_i
   < \infty$.
Indeed, it holds that $e_i = u_i - p$ for any
$i \in \Lambda({\rm GS}^2(E_1,p))$, where
$u_i \in {\rm GS}^2(E_1,p) \cap B_r(p)$.

Let $y$ be an arbitrary element of ${\rm GS}^2(E_2,q)$.
Since
$T_{-q} \circ F_2 \circ T_p : {\rm GS}^2(E_1,p) - p \to
 {\rm GS}^2(E_2,q) - q$ is a surjective $d_a$-isometry,
there exists an element $x \in {\rm GS}^2(E_1,p)$ with
$y - q = (T_{-q} \circ F_2 \circ T_p)(x - p)$, where according
to Theorem \ref{thm:3.7}, there are some real numbers
$\alpha_i$ such that
\begin{eqnarray} \label{eq:2021-11-30}
x - p =   \sum_{i \in \Lambda({\rm GS}^2(E_1,p))}
          \alpha_i e_i
      =   \sum_{i \in \Lambda({\rm GS}^2(E_1,p))}
          \alpha_i (u_i - p)
      \in {\rm GS}^3(E_1,p) - p.
\end{eqnarray}

On the other hand, due to Definition \ref{def:6.1} $(ii)$ and
(\ref{eq:2021-11-30}), it holds that
\begin{eqnarray} \label{eq:2021-12-1}
\big( T_{-q} \circ F_3 \circ T_p \big) (x - p)
= \sum_{i \in \Lambda({\rm GS}^2(E_1,p))}
  \alpha_i \big( T_{-q} \circ F_2 \circ T_p \big) (u_i - p).
\end{eqnarray}
Furthermore, since $e_i = u_i - p$ for all
$i \in \Lambda({\rm GS}^2(E_1,p))$, it follows from Theorem
\ref{thm:setyol}, (\ref{eq:2021-11-30}), and
(\ref{eq:2021-12-1}) that
\begin{eqnarray} \label{eq:2021-11-30-1}
\begin{split}
y - q
& = \big( T_{-q} \circ F_2 \circ T_p \big)\!
    \left( \sum_{i \in \Lambda({\rm GS}^2(E_1,p))} \alpha_i e_i
    \right) \\
& = \sum_{i \in \Lambda({\rm GS}^2(E_1,p))}
    \alpha_i \big( T_{-q} \circ F_2 \circ T_p \big) (e_i).
\end{split}
\end{eqnarray}

Due to Lemma \ref{lem:6.1},
$\{ \frac{1}{a_i} (T_{-q} \circ F_2 \circ T_p)(e_i) \}
 _{i \in \Lambda({\rm GS}^2(E_1,p))}$ is an orthonormal
sequence in ${\rm GS}^2(E_2,q) - q$.
Assume that
$\langle y - q, \frac{1}{a_i} (T_{-q} \circ F_2 \circ T_p)(e_i)
 \rangle_a = 0$ for all $i \in \Lambda({\rm GS}^2(E_1,p))$.
It then follows from (\ref{eq:2021-11-30-1}) that
$\alpha_i = 0$ for each $i \in \Lambda({\rm GS}^2(E_1,p))$.
Then, by (\ref{eq:2021-11-30-1}) again, we have
\begin{eqnarray*}
y - q = \big( T_{-q} \circ F_2 \circ T_p \big) (0) = 0,
\end{eqnarray*}
which implies that
$\{ \frac{1}{a_i} (T_{-q} \circ F_2 \circ T_p)(e_i) \}
 _{i \in \Lambda({\rm GS}^2(E_1,p))}$ is a complete orthonormal
sequence in ${\rm GS}^2(E_2,q) - q$.
\hfill$\Box$
\vspace*{5mm}

Assume that there exists an index $i_0$ of
${\rm GS}^2(\tilde{J},p)$ but it is not an index of $\tilde{J}$.
According to Theorem \ref{thm:3.7}, we have
\begin{eqnarray*}
p + e_{i_0} \in {\rm GS}^2(\tilde{J},p).
\end{eqnarray*}
In view of Definitions \ref{def:gls} and \ref{def:6.1} $(i)$,
there exist an $m_0 \in \mathbb{N}$, some
$x_{ij} \in {\rm GS}(\tilde{J},p) \cap B_r(p)$, and some real
numbers $\alpha_{ij}$ such that
\begin{eqnarray*}
p + e_{i_0}
=   p +
    \sum_{i=1}^{m_0} \sum_{j=1}^\infty \alpha_{ij} (x_{ij} - p)
\in M_a.
\end{eqnarray*}

Since each $x_{ij}$ belongs to ${\rm GS}(\tilde{J},p)$, it
follows from Theorem \ref{thm:set9} that
\begin{eqnarray*}
x_{ij} = p + \sum_{k \in \Lambda(\tilde{J})} \beta_{ijk} e_k
\end{eqnarray*}
for all $i, j$ and for some real numbers $\beta_{ijk}$.
Hence, we get
\begin{eqnarray*}
p + e_{i_0}
= p +
  \sum_{i=1}^{m_0} \sum_{j=1}^\infty \alpha_{ij}
  \sum_{k \in \Lambda(\tilde{J})} \beta_{ijk} e_k
\end{eqnarray*}
or
\begin{eqnarray*}
\frac{1}{a_{i_0}} e_{i_0}
= \sum_{i=1}^{m_0} \sum_{j=1}^\infty
    \sum_{k \in \Lambda(\tilde{J})}
    \frac{a_k}{a_{i_0}} \alpha_{ij} \beta_{ijk}
    \frac{1}{a_k} e_k.
\end{eqnarray*}

We remind that $i_0 \not\in \Lambda(\tilde{J})$ and
$\{ \frac{1}{a_i} e_i \}_{i \in \mathbb{N}}$ is an orthonormal
sequence in $M_a$.
And we obtain
\begin{eqnarray*}
1 = \left\langle \frac{1}{a_{i_0}} e_{i_0},
                 \frac{1}{a_{i_0}} e_{i_0}
    \right\rangle_a
  = \sum_{i=1}^{m_0} \sum_{j=1}^\infty
    \sum_{k \in \Lambda(\tilde{J})}
    \frac{a_k}{a_{i_0}} \alpha_{ij} \beta_{ijk}
    \left\langle \frac{1}{a_k} e_k,
                 \frac{1}{a_{i_0}} e_{i_0}
    \right\rangle_a
  = 0,
\end{eqnarray*}
which leads to the contradiction.
Therefore, we conclude that
$\Lambda({\rm GS}^2(\tilde{J},p)) = \Lambda(\tilde{J})$.

This fact will be important to prove the following theorem.
It is especially important for deriving the second expression
in (\ref{eq:8.4}).


\begin{lem} \label{lem:8.2}
Assume that $E_1$ is a bounded subset of $I^\omega$ that
contains at least two elements and is $d_a$-isometric to a
subset $E_2$ of $I^\omega$ via a surjective $d_a$-isometry
$f : E_1 \to E_2$.
Moreover, assume that $p$ and $q$ are elements of $E_1$ and
$E_2$ that satisfy $q = f(p)$.
Assume that $\Lambda$ is a set with
$\Lambda({\rm GS}^2(E_1,p)) \subset \Lambda \subset \mathbb{N}$
and $\{ \frac{1}{a_i} e'_i \}_{i \in \Lambda}$ is the
orthonormal sequence given in Definition \ref{def:3.5}.
Then $\{ \frac{1}{a_i} e'_i \}_{i \in \Lambda}$ is a complete
orthonormal sequence in $G_2({\rm GS}^2(\tilde{J},p)) - q$.
\end{lem}

\noindent
\emph{Proof}.
If we examine the structure of $\tilde{J}$ given in Definition
\ref{def:3.5}, we see that $\Lambda = \Lambda(\tilde{J})$.
It then follows from Theorem \ref{thm:set9} that
$p + e_i \in     {\rm GS}(\tilde{J},p)
         \subset {\rm GS}^2(\tilde{J},p)$ for any
$i \in \Lambda(\tilde{J})$.
We assume that $x \in {\rm GS}^2(\tilde{J},p)$.
Then, according to Theorem \ref{thm:3.7}, there are some real
numbers $\alpha_i$ such that
\begin{eqnarray} \label{eq:8.4}
x - p =   \sum_{i \in \Lambda({\rm GS}^2(\tilde{J},p))}
          \alpha_i e_i
      =   \sum_{i \in \Lambda(\tilde{J})} \alpha_i e_i.
\end{eqnarray}

For any $i \in \Lambda(\tilde{J}) = \Lambda$, it follows from
(\ref{eq:2021-10-21-1}) that
\begin{eqnarray*}
\big( T_{-q} \circ G_2 \circ T_p \big) (e_i)
= \sum_{j \in \Lambda(\tilde{J})}
  \left\langle e_i, \frac{1}{a_j} e_j \right\rangle_a
  \frac{1}{a_j} e'_j
= e'_i,
\end{eqnarray*}
which implies that the real Hilbert space
$G_2({\rm GS}^2(\tilde{J},p)) - q$ contains the orthonormal
sequence $\{ \frac{1}{a_i} e'_i \}_{i \in \Lambda}$.

Let $y$ be an arbitrary element of
$G_2({\rm GS}^2(\tilde{J},p))$.
Since
$T_{-q} \circ G_2 \circ T_p : {\rm GS}^2(\tilde{J},p) - p \to
 G_2({\rm GS}^2(\tilde{J},p)) - q$ is a surjective
$d_a$-isometry, there exists an element
$x \in {\rm GS}^2(\tilde{J},p)$ with
$y - q = (T_{-q} \circ G_2 \circ T_p)(x - p)$, where $x - p$
is expressed by (\ref{eq:8.4}).
Then, by (\ref{eq:2021-10-21-1}) and (\ref{eq:8.4}), we get
\begin{eqnarray} \label{eq:8.5}
\begin{split}
y - q
& = \big( T_{-q} \circ G_2 \circ T_p \big)\!
    \left( \sum_{i \in \Lambda(\tilde{J})} \alpha_i e_i
    \right) \\
& = \sum_{j \in \Lambda(\tilde{J})}
    \left\langle \sum_{i \in \Lambda(\tilde{J})} \alpha_i e_i,
                 \frac{1}{a_j} e_j
    \right\rangle_a \frac{1}{a_j} e'_j \\
& = \sum_{j \in \Lambda(\tilde{J})} \alpha_j e'_j.
\end{split}
\end{eqnarray}

Assume that
$\langle y - q, \frac{1}{a_i} e'_i \rangle_a = 0$ for all
$i \in \Lambda(\tilde{J})$.
It then follows from (\ref{eq:8.5}) that $\alpha_i = 0$ for
each $i \in \Lambda(\tilde{J})$.
Then, by (\ref{eq:8.5}) again, we have
\begin{eqnarray*}
y - q = \big( T_{-q} \circ G_2 \circ T_p \big) (0) = 0,
\end{eqnarray*}
which implies that
$\{ \frac{1}{a_i} e'_i \}_{i \in \Lambda}$ is a complete
orthonormal sequence in $G_2({\rm GS}^2(\tilde{J},p)) - q$.
\hfill$\Box$
\vspace*{5mm}

The following theorem states that the domain of local
$d_a$-isometry can be extended to any real Hilbert space
including the domain of local $d_a$-isometry.


\begin{thm} \label{thm:3.17}
Let $E_1$ be a bounded subset of $I^\omega$ that contains at
least two elements.
Suppose $E_1$ is $d_a$-isometric to a subset $E_2$ of
$I^\omega$ via a surjective $d_a$-isometry $f : E_1 \to E_2$.
Let $p$ and $q$ be elements of $E_1$ and $E_2$ satisfying
$q = f(p)$.
Assume that
$\mathbb{N} \!\setminus\! \Lambda({\rm GS}^2(E_1,p))$ is a
nonempty set and $\Lambda$ is a set with
$\Lambda({\rm GS}^2(E_1,p)) \subset \Lambda \subset \mathbb{N}$.
Suppose $\{ \frac{1}{a_i} e'_i \}_{i \in \Lambda}$ is an
orthonormal sequence and
$e'_i = (T_{-q} \circ F_2 \circ T_p)(e_i)$ for each
$i \in \Lambda({\rm GS}^2(E_1,p))$.
Let $\tilde{J}$ be either a basic cylinder or a basic interval
defined in Definition \ref{def:3.5}.
Then the function $G_2 : {\rm GS}^2(\tilde{J},p) \to M_a$ is a
$d_a$-isometry and the function
$T_{-q} \circ G_2 \circ T_p : {\rm GS}^2(\tilde{J},p) - p \to M_a$
is linear.
In particular, $G_2$ is an extension of $F_2$.
\end{thm}

\noindent
\emph{Proof}.
$(a)$
First, we assert that the function
$T_{-q} \circ G_2 \circ T_p : {\rm GS}^2(\tilde{J},p) - p \to M_a$
preserves the inner product.
Assume that $u$ and $v$ are arbitrary elements of
${\rm GS}^2(\tilde{J},p)$.
Since $\Lambda = \Lambda(\tilde{J})$, it follows from
(\ref{eq:2021-10-21-2}), (\ref{eq:2021-10-21-1}), and the
orthonormality of $\{ \frac{1}{a_i} e_i \}$ and
$\{ \frac{1}{a_i} e'_i \}$ that
\begin{eqnarray*}
\begin{split}
& \big\langle (T_{-q} \circ G_2 \circ T_p)(u - p),\,
              (T_{-q} \circ G_2 \circ T_p)(v - p)
  \big\rangle_a \\
& \hspace*{3mm}
  = \left\langle \sum_{i \in \Lambda}
                 \left\langle u - p,\, \frac{1}{a_i} e_i
                 \right\rangle_a
                 \frac{1}{a_i} e'_i,\,
                 \sum_{j \in \Lambda}
                 \left\langle v - p,\, \frac{1}{a_j} e_j
                 \right\rangle_a
                 \frac{1}{a_j} e'_j
    \right\rangle_a \\
& \hspace*{3mm}
  = \sum_{i \in \Lambda}
    \left\langle u - p,\, \frac{1}{a_i} e_i \right\rangle_a
    \sum_{j \in \Lambda}
    \left\langle v - p,\, \frac{1}{a_j} e_j \right\rangle_a
    \left\langle \frac{1}{a_i} e'_i,\, \frac{1}{a_j} e'_j
    \right\rangle_a \\
& \hspace*{3mm}
  = \sum_{i \in \Lambda}
    \left\langle u - p,\, \frac{1}{a_i} e_i \right\rangle_a
    \sum_{j \in \Lambda}
    \left\langle v - p,\, \frac{1}{a_j} e_j \right\rangle_a
    \left\langle \frac{1}{a_i} e_i,\, \frac{1}{a_j} e_j
    \right\rangle_a \\
& \hspace*{3mm}
  = \left\langle \sum_{i \in \Lambda}
                 \left\langle u - p,\, \frac{1}{a_i} e_i
                 \right\rangle_a
                 \frac{1}{a_i} e_i,\,
                 \sum_{j \in \Lambda}
                 \left\langle v - p,\, \frac{1}{a_j} e_j
                 \right\rangle_a
                 \frac{1}{a_j} e_j
    \right\rangle_a \\
& \hspace*{3mm}
  = \langle u - p, v - p \rangle_a
\end{split}
\end{eqnarray*}
for all $u, v \in {\rm GS}^2(\tilde{J},p)$, \textit{i.e.},
$T_{-q} \circ G_2 \circ T_p$ preserves the inner product.

$(b)$
We assert that $G_2$ is a $d_a$-isometry.
Let $u$ and $v$ be arbitrary elements of
${\rm GS}^2(\tilde{J},p)$.
Since $T_{-q} \circ G_2 \circ T_p$ preserves the inner product
by $(a)$, we have
\begin{eqnarray*}
\begin{split}
& d_a \big( G_2(u), G_2(v) \big)^2 \\
& \hspace*{3mm}
  = \big\| ( T_{-q} \circ G_2 \circ T_p )(u-p) -
           ( T_{-q} \circ G_2 \circ T_p )(v-p)
    \big\|_a^2 \\
& \hspace*{3mm}
  = \big\langle ( T_{-q} \circ G_2 \circ T_p )(u-p) -
                ( T_{-q} \circ G_2 \circ T_p )(v-p), \\
& \hspace*{3mm}
    ~~~~~~      ( T_{-q} \circ G_2 \circ T_p )(u-p) -
                ( T_{-q} \circ G_2 \circ T_p )(v-p)
    \big\rangle_a \\
& \hspace*{3mm}
  = \langle u-p, u-p \rangle_a - \langle u-p, v-p \rangle_a
    - \langle v-p, u-p \rangle_a + \langle v-p, v-p \rangle_a \\
& \hspace*{3mm}
  = \big\langle (u-p) - (v-p), (u-p) - (v-p) \big\rangle_a \\
& \hspace*{3mm}
  = \| (u-p) - (v-p) \|_a^2 \\
& \hspace*{3mm}
  = \| u - v \|_a^2 \\
& \hspace*{3mm}
  = d_a(u, v)^2
\end{split}
\end{eqnarray*}
for all $u, v \in {\rm GS}^2(\tilde{J},p)$, \textit{i.e.},
$G_2 : {\rm GS}^2(\tilde{J},p) \to M_a$ is a $d_a$-isometry.

$(c)$
Now we assert that the function
$T_{-q} \circ G_2 \circ T_p : {\rm GS}^2(\tilde{J},p) - p \to M_a$
is linear.
Assume that $u$ and $v$ are arbitrary elements of
${\rm GS}^2(\tilde{J},p)$ and $\alpha$ and $\beta$ are real
numbers.
Since ${\rm GS}^2(\tilde{J},p) - p$ is a real vector space, it
holds that 
$\alpha (u - p) + \beta (v - p) \in {\rm GS}^2(\tilde{J},p) - p$.
Thus, $\alpha (u - p) + \beta (v - p) = w - p$ for some
$w \in {\rm GS}^2(\tilde{J},p)$.
Hence, referring to the changes presented in the table below
and following $(d)$ of the proof of Theorem \ref{thm:new3}, we
can easily prove that $T_{-q} \circ G_2 \circ T_p$ is linear.

\begin{center}
\begin{tabular}{cccc}
\hline\hline
Theorem \ref{thm:new3}: & ${\rm GS}(E_1,p)$         & $F$
& (\ref{eq:2018-9-28}) \\
\hline
Here:                   & ${\rm GS}^2(\tilde{J},p)$ & $G_2$
& $(a)$  \\
\hline\hline
\end{tabular}
\end{center}

$(d)$
Finally, we assert that $G_2$ is an extension of $F_2$.
Let $\hat{J}$ be either a basic cylinder or a basic interval
defined by
\begin{eqnarray*}
\hat{J} = \prod_{i=1}^\infty \hat{J}_i, ~~~\mbox{where}~~~
\hat{J}_i
= \left\{\begin{array}{ll}
            [0, 1]
            & (\mbox{for}~i \in \Lambda({\rm GS}^2(E_1,p))),
            \vspace*{1mm}\\
            \{ p_i \}
            & (\mbox{for}~i \not\in \Lambda({\rm GS}^2(E_1,p))).
         \end{array}
  \right.
\end{eqnarray*}
We see that $p = (p_1, p_2, \ldots) \in \hat{J} \cap E_1$ and
$\Lambda(\hat{J}) = \Lambda({\rm GS}^2(E_1,p))$.
Due to Theorem \ref{thm:setten}, we conclude that
$\hat{J} \subset {\rm GS}(\hat{J},p) \subset {\rm GS}^2(E_1,p)$.
So we have
\begin{eqnarray*}
\hat{J} \cap B_r(p) \subset {\rm GS}(\hat{J},p) \cap B_r(p)
\subset {\rm GS}^2(E_1,p) \cap B_r(p)
\end{eqnarray*}
for some real number $r > 0$ and hence, we further have
\begin{eqnarray*}
{\rm GS}(\hat{J},p) \subset {\rm GS}^2(\hat{J},p)
\subset {\rm GS}^3(E_1,p) = {\rm GS}^2(E_1,p).
\end{eqnarray*}

On the other hand, by Theorems \ref{thm:set9} and \ref{thm:3.7},
we have
\begin{eqnarray*}
\begin{split}
{\rm GS}^2(E_1,p)
& = \Bigg\{ p + \sum_{i \in \Lambda({\rm GS}^2(E_1,p))}
            \alpha_i e_i \in M_a :
            \alpha_i \in \mathbb{R} ~\mbox{for all}~
            i \in \Lambda({\rm GS}^2(E_1,p))
    \Bigg\} \\
& = \Bigg\{ p + \sum_{i \in \Lambda(\hat{J})}
            \alpha_i e_i \in M_a :
            \alpha_i \in \mathbb{R} ~\mbox{for all}~
            i \in \Lambda(\hat{J})
    \Bigg\} \\
& = {\rm GS}(\hat{J},p) \\
& \subset {\rm GS}^2(\hat{J},p) \\
& \subset {\rm GS}^2(E_1,p),
\end{split}
\end{eqnarray*}
which implies that
${\rm GS}(\hat{J},p) = {\rm GS}^2(\hat{J},p) = {\rm GS}^2(E_1,p)$.

Let $u$ be an arbitrary element of
${\rm GS}^2(E_1,p) = {\rm GS}^2(\hat{J},p)$.
Then by (\ref{eq:2021-10-21-2}) with $\hat{J}$ instead of $J$,
we have
\begin{eqnarray} \label{eq:2021-11-8}
u - p = \sum_{i \in \Lambda(\hat{J})}
        \left\langle u - p,\, \frac{1}{a_i} e_i \right\rangle_a
        \frac{1}{a_i} e_i
\end{eqnarray}
and since $T_{-q} \circ F_2 \circ T_p$ is linear and continuous,
we apply Lemma \ref{lem:doolset} with $\hat{J}$ in place of $E$,
and we use (\ref{eq:2021-10-21-1}), (\ref{eq:2021-11-8}), and
the facts
${\rm GS}^2(E_1,p) = {\rm GS}^2(\hat{J},p) = {\rm GS}(\hat{J},p)$
and $\Lambda(\hat{J}) = \Lambda({\rm GS}^2(E_1,p))$ to have
\begin{eqnarray*}
\begin{split}
(T_{-q} \circ G_2 \circ T_p)(u - p)
& = \sum_{i \in \Lambda(\tilde{J})}
    \left\langle u - p,\, \frac{1}{a_i} e_i \right\rangle_a
    \frac{1}{a_i} e'_i \\
& = \sum_{i \in \Lambda(\hat{J})}
    \left\langle u - p,\, \frac{1}{a_i} e_i \right\rangle_a
    \frac{1}{a_i} e'_i \\
& = \sum_{i \in \Lambda(\hat{J})}
    \left\langle u - p,\, \frac{1}{a_i} e_i \right\rangle_a
    \frac{1}{a_i} (T_{-q} \circ F_2 \circ T_p)(e_i) \\
& = \lim_{n \to \infty}
    \sum_{i \in \Lambda_n(\hat{J})}
    \left\langle u - p,\, \frac{1}{a_i} e_i \right\rangle_a
    \frac{1}{a_i} (T_{-q} \circ F_2 \circ T_p)(e_i) \\
& = \lim_{n \to \infty}
    (T_{-q} \circ F_2 \circ T_p)\!
    \left(
    \sum_{i \in \Lambda_n(\hat{J})}
    \left\langle u - p,\, \frac{1}{a_i} e_i \right\rangle_a
    \frac{1}{a_i} e_i
    \right) \\
& = (T_{-q} \circ F_2 \circ T_p)\!
    \left(
    \sum_{i \in \Lambda(\hat{J})}
    \left\langle u - p,\, \frac{1}{a_i} e_i \right\rangle_a
    \frac{1}{a_i} e_i
    \right) \\
& = (T_{-q} \circ F_2 \circ T_p)(u - p),
\end{split}
\end{eqnarray*}
where we set
$\Lambda_n(\hat{J}) = \{ i \in \Lambda(\hat{J}) : i < n \}$
for every $n \in \mathbb{N}$.

Therefore, it follows that
$(T_{-q} \circ G_2 \circ T_p)(u - p)
 = (T_{-q} \circ F_2 \circ T_p)(u - p)$ for all
$u \in {\rm GS}^2(E_1,p)$, \textit{i.e.}, $G_2(u) = F_2(u)$
for all $u \in {\rm GS}^2(E_1,p)$.
In other words, $G_2$ is an extension of $F_2$.
\hfill$\Box$


\section{Applications}
\label{sec:9}

Given an integer $n > 0$, let $a = \{ a_i \}$ be a sequence of
positive real numbers that satisfies
\begin{eqnarray*}
a_1 = \cdots = a_n = 1 ~~~\mbox{and}~~~
\sum_{i=n+1}^\infty a_i^2 < \infty.
\end{eqnarray*}
We know that the $n$ dimensional real vector space
$\mathbb{R}^n$ can be identified with a subspace of $M_a$.
More precisely, it holds that
$\mathbb{R}^n \simeq
 \big\{ (x_1, \ldots, x_n, 0, 0, \ldots) : x_i \in \mathbb{R}
        ~\mbox{for}~ 1 \leq i \leq n
 \big\}$.

We can define an inner product $\langle \cdot, \cdot \rangle_a$
for $\mathbb{R}^n$ by
\begin{eqnarray*}
\langle x, y \rangle_a = \sum_{i=1}^\infty a_i^2 x_i y_i
                       = \sum_{i=1}^n x_i y_i
\end{eqnarray*}
for all $x, y \in \mathbb{R}^n$, with which
$(\mathbb{R}^n, \langle \cdot, \cdot \rangle_a)$ becomes a real
Hilbert space.
This inner product induces the Euclidean norm in the natural way
\begin{eqnarray*}
\| x \|_a = \sqrt{\langle x, x \rangle_a}
          = \sqrt{\sum_{i=1}^n x_i^2}
\end{eqnarray*}
for all $x \in \mathbb{R}^n$, and
$(\mathbb{R}^n, \| \cdot \|_a)$ becomes a real Banach space.

If we replace $M_a$, $\mathbb{N}$ and $I^\omega$ with
$\mathbb{R}^n$, $\{ 1, 2, \ldots, n \}$ and $[0, 1]^n$,
respectively, in the definitions and theorems in the previous
sections, it is not difficult to see that they also hold for
the $n$ dimensional Euclidean space $\mathbb{R}^n$.

The following theorem is the finite dimensional real Hilbert
space version of Theorem \ref{thm:3.17}, the main theorem of
this paper.
More specifically, Theorem \ref{thm:3.17}, which applies to
infinite dimensional real Hilbert spaces, is also applicable
to finite dimensional real Hilbert spaces.


\begin{thm} \label{thm:9.1}
Let $E_1$ be a bounded subset of\/ $[0, 1]^n$ that contains at
least two elements.
Suppose $E_1$ is $d_a$-isometric to a subset $E_2$ of\/
$[0, 1]^n$ via a surjective $d_a$-isometry $f : E_1 \to E_2$.
Let $p$ and $q$ be elements of $E_1$ and $E_2$ satisfying
$q = f(p)$.
Assume that $\Lambda$ is a finite set with
$\Lambda({\rm GS}^2(E_1,p)) \subset \Lambda \subset
 \{ 1, 2, \ldots, n \}$.
Suppose $\{ e'_i \}_{i \in \Lambda}$ are orthonormal vectors
and $e'_i = (T_{-q} \circ F_2 \circ T_p)(e_i)$ for each
$i \in \Lambda({\rm GS}^2(E_1,p))$.
Let $\tilde{J}$ be a basic interval defined in Definition
\ref{def:3.5}.
Then the function
$G_2 : {\rm GS}^2(\tilde{J},p) \to \mathbb{R}^n$ is a
$d_a$-isometry and the function
$T_{-q} \circ G_2 \circ T_p : {\rm GS}^2(\tilde{J},p) - p \to
 \mathbb{R}^n$ is linear.
In particular, $G_2$ is an extension of $F_2$.
\end{thm}


\section{Discussion}
\label{sec:10}

In view of Lemma \ref{lem:3.8} $(iii)$,
${\rm GS}^2(\tilde{J},p) - p$ is closed in $M_a$.
Hence, ${\rm GS}^2(\tilde{J},p) - p$ is a closed subspace of
the real Hilbert space $M_a$, \textit{i.e.},
${\rm GS}^2(\tilde{J},p) - p$ is itself a real Hilbert space.
In addition, since the $d_a$-isometry
$T_{-q} \circ G_2 \circ T_p : {\rm GS}^2(\tilde{J},p) - p \to
 G_2({\rm GS}^2(\tilde{J},p)) - q$ is a homeomorphism and
$(T_{-q} \circ G_2 \circ T_p)({\rm GS}^2(\tilde{J},p) - p)
 = G_2({\rm GS}^2(\tilde{J},p)) - q$,
$G_2({\rm GS}^2(\tilde{J},p)) - q$ is a closed subspace of the
real Hilbert space $M_a$, \textit{i.e.}, it is also a real
Hilbert space.

In view of Lemma \ref{lem:8.2}, we notice that
$\{ \frac{1}{a_i} e'_i \}_{i \in \Lambda}$ is a complete
orthonormal sequence in $G_2({\rm GS}^2(\tilde{J},p)) - q$
and, since ${\rm GS}^2(E_2,q) - q$ is a subspace of
$G_2({\rm GS}^2(\tilde{J},p)) - q$, it follows from
\cite[Theorems 3.9.1 and 3.9.4]{debnath} that
\begin{eqnarray*}
\begin{split}
& G_2 \big( {\rm GS}^2(\tilde{J},p) \big) - q \\
& \hspace*{3mm}
  = ({\rm GS}^2(E_2,q) - q) \oplus
    \big( ({\rm GS}^2(E_2,q) - q)^\perp \cap
          (G_2({\rm GS}^2(\tilde{J},p)) - q)
    \big).
\end{split}
\end{eqnarray*}
More precisely, the real Hilbert space
$G_2 \big( {\rm GS}^2(\tilde{J},p) \big) - q$ is the direct
sum of two real Hilbert spaces ${\rm GS}^2(E_2,q) - q$ and
$({\rm GS}^2(E_2,q) - q)^\perp \cap
 (G_2({\rm GS}^2(\tilde{J},p)) - q)$.

It is well known that every Hilbert space has a complete
orthonormal sequence and
the union of a complete orthonormal sequence in the Hilbert
subspace ${\rm GS}^2(E_2,q) - q$ and a complete orthonormal
sequence in the orthogonal complement of ${\rm GS}^2(E_2,q) - q$
(as a subspace of $G_2({\rm GS}^2(\tilde{J},p)) - q$) is a
complete orthonormal sequence in the Hilbert space
$G_2({\rm GS}^2(\tilde{J},p)) - q$.
Hence, in view of Lemmas \ref{lem:8.1} and \ref{lem:8.2},
\begin{eqnarray*}
\left\{ \frac{1}{a_i} e'_i \right\}_{i \in \Lambda}
= \left\{ \frac{1}{a_i} (T_{-q} \circ F_2 \circ T_p)(e_i)
  \right\}_{i \in \Lambda({\rm GS}^2(E_1,p))} \cup
  \left\{ \frac{1}{a_i} e'_i
  \right\}_{i \in \Lambda \setminus \Lambda({\rm GS}^2(E_1,p))}
\end{eqnarray*}
is a complete orthonormal sequence in the real Hilbert space
$G_2({\rm GS}^2(\tilde{J},p)) - q$.

From the above facts, together with Theorem \ref{thm:3.17}, we
know that the domain of local $d_a$-isometry can be extended
to the real Hilbert space containing that domain.


\section{Conclusion}
\label{sec:11}

In view of Theorem \ref{thm:3.17}, the domain of local
$d_a$-isometry can be extended to the real Hilbert space
containing that domain.
The domain of local $d_a$-isometry does not need to be a convex
body or an open set required by \cite{manki}, it just needs to
be bounded and contain at least two elements.
Therefore, the coverage of our result is wider than that of
previous results.
This is the biggest advantage of this paper compared to the
previous results.
\vspace*{5mm}


\section*{Acknowledgments}
This work was supported by the National Research Foundation of
Korea (NRF) grant funded by the Korea government (MSIT)
(No. 2020R1F1A1A01049560).

\end{document}